\documentclass[11pt]{article}%
\usepackage{amsmath}
\usepackage{graphicx}
\usepackage{amsfonts}
\usepackage{amssymb}
\usepackage{a4wide}%
\usepackage[colorlinks=true,breaklinks=true,bookmarks=true,urlcolor=blue,
     citecolor=blue,linkcolor=blue,bookmarksopen=false,draft=false]{hyperref}
\providecommand{\U}[1]{\protect\rule{.1in}{.1in}}
\newtheorem{theorem}{Theorem}

\newtheorem{corollary}[theorem]{Corollary}

\newtheorem{example}{Example}

\newtheorem{lemma}[theorem]{Lemma}

\newtheorem{proposition}[theorem]{Proposition}
\newtheorem{remark}{Remark}

\begin{document}

\title{Relations between the convexity of a set and the differentiability of its
support function}
\author{C. Z\u{a}linescu\thanks{University Alexandru Ioan Cuza, Faculty of
Mathematics,  Ia\c{s}i, Romania, email: \texttt{zalinesc@uaic.ro}.}}
\date{}
\maketitle

\begin{abstract}
It is known that, in finite dimensions, the support function of a
compact convex set with non empty interior is differentiable
excepting the origin if and only if the set is strictly convex. In
this paper we realize a thorough study of the relations between the
differentiability of the support function on the interior of its
domain and the convexity of the set, mainly for unbounded sets. Then
we revisit some results related to the differentiability of the cost
function associated to a production function.
\end{abstract}

\section{Introduction}

The celebrated Shephard lemma, which is considered to be ``a major result in
microeconomics having applications in the theory of the firm and in consumer
choice'' (see http://en.wikipedia.\allowbreak org/wiki/Shephard's\_lemma) is
related to the differentiability of the cost function in economics. The cost
function is defined by%
\[
g:\mathbb{R}_{++}^{p}\rightarrow\mathbb{R},\quad g(x):=\inf\left\{
\left\langle x,a\right\rangle \mid a\in A\right\}  ,
\]
where $A$ is a nonempty subset of $\mathbb{R}_{+}^{p}$. More precisely,
$A=\{u\in\mathbb{R}_{+}^{p}\mid f(u)\geq y\}$, where $f:\mathbb{R}_{+}%
^{p}\rightarrow\mathbb{R}_{+}$ is a production function.

Clearly, the function $g$ above is strongly related to the support function
which is defined as
\[
\sigma_{A}:X^{\ast}\rightarrow\overline{\mathbb{R}},\quad\sigma_{A}(x^{\ast
}):=\sup\left\{  \left\langle x^{\ast},u\right\rangle \mid u\in A\right\}  ,
\]
where $X$ is a (finite dimensional) real normed space whose topological dual
is denoted by $X^{\ast}$, and $A\subset X$ is a nonempty set.

Because $g(x^{\ast})=-\sigma_{A}(-x^{\ast})$ for $x^{\ast}\in\mathbb{R}^{p}$
(where $X=\mathbb{R}^{p}$ is endowed with the usual Euclidean norm), any
property of the support function $\sigma_{A}$ can be translated into a
corresponding property of the cost function $g$.

In the economics literature one can find several results related to
Shephard's lemma and to the differentiability of the cost function;
see \cite{FarPri:86}, \cite{FarPriSam:90}, \cite{Fuc:95},
\cite{Fuc:97}, \cite{Kim:93}, \cite{Sai:83}, \cite{Sak:73}. Our aim
is to study the connection between the differentiability of the
support function $\sigma _A$ and the convexity of $A$, and  to
revisit some results related to Shephard's lemma.

The rest of the paper is organized as follows. In Section 2 we
present preliminary notions and results. In Section 3 we recall
several results from \cite{Zal:12} concerning the differentiability
of $\sigma _A$ for $A$ convex. Section 4 contains the main results
of this paper. We present several conditions on the set $A$ which
imply the differentiability of $\sigma _A$ (on certain sets) and
conditions on $A$ under which the differentiability of $\sigma _A$
on $\operatorname*{int}(\operatorname*{dom}\sigma _A)$ implies the
convexity of $A$. Moreover, we associate to a set $A$ satisfying
condition (H) a function $F_A$ and establish relationships between
properties of the set $A$ and properties of $F_A$. Then, in Section
5, we apply the results in Section 4 to the problem of the
differentiability of the cost function and discuss several results
on Shephard's lemma from the economics literature.

\section{Preliminaries}

In the following we assume that $X$ is a nontrivial real finite
dimensional normed space whose dual is denoted by $X^{\ast}$. We
identify $(X^{\ast })^{\ast}$ with $X$. However, the reader can take
$X$ an Euclidean space and identify $X^{\ast}$ with $X$. For
$A\subset X$ we denote by $\operatorname*{aff}A$,
$\operatorname*{lin}_{0}A$, $\operatorname*{cl}A$,
$\operatorname*{int}A$, $\operatorname*{rint}A$,
$\operatorname*{bd}A$, $\operatorname*{rbd}A,$
$\operatorname*{conv}A$, $\overline {\operatorname*{conv}}A$ the
affine hull of $A$, the linear space parallel to
$\operatorname*{aff}A$, the closure of $A$, the interior of $A$, the
relative interior of $A$ (that is the interior of $A$ w.r.t.\
$\operatorname*{aff}A$),
the boundary of $A$ (hence $\operatorname*{bd}A=\operatorname*{cl}%
A\setminus\operatorname*{int}A)$, the relative boundary of $A$ (hence
$\operatorname*{rbd}A=\operatorname*{cl}A\setminus\operatorname*{rint}A)$, the
convex hull of $A$ and $\operatorname*{cl}(\operatorname*{conv}A)$,
respectively. When $A\subset B\subset X$, we write $\operatorname*{int}_{B}A$
and $\operatorname*{bd}_{B}A$ for the interior and boundary of $A$ as subset
of $B$ endowed with the induced topology. Recall that the recession cone of
$A\neq\emptyset$ is the set
\[
A_{\infty}:=\left\{  u\in X\mid\exists(t_{n})\subset(0,\infty),\ t_{n}%
\rightarrow0,\ \exists(a_{n})\subset A:t_{n}a_{n}\rightarrow u\right\}  .
\]
Clearly, $A_{\infty}$ is a (closed) cone; in particular, $0\in A_{\infty}$. We
have that $A_{\infty}=\{0\}$ if and only if $A$ is bounded. When $A$ is closed
and convex we have that $A_{\infty}=\cap_{t>0}t(A-a)$, where $a\in A$. In this
case $A_{\infty}$ is also convex; moreover, $A_{\infty}$ is pointed, that is,
$A_{\infty}\cap(-A_{\infty})=\{0\}$, if and only if $A$ does not contain any line.

For the the set $A\subset X$ we set
\[
A^{+}:=\{x^{\ast}\in X^{\ast}\mid\left\langle x,x^{\ast}\right\rangle
\geq0\ \forall x\in A\},\quad A^{-}:=-A^{+},\quad A^{\perp}:=A^{+}\cap A^{-},
\]
where $\left\langle x,x^{\ast}\right\rangle :=x^{\ast}(x)$ for $x\in X$ and
$x^{\ast}\in X^{\ast}$.

For $P\subset X$ a closed convex cone we set
\[
P^{\#}:=\{x^{\ast}\in X^{\ast}\mid\left\langle x,x^{\ast}\right\rangle
>0\ \forall x\in P\setminus\{0\}\}.
\]
It is known that $P^{\#}\neq\emptyset$ if and only if $P$ is pointed (that is
$P\cap(-P)=\{0\}$); note that $\{0\}^{\#}=X^{\ast}$. Furthermore,
$P^{\#}=\operatorname*{int}P^{+}$. Moreover, for $x,x^{\prime}\in X$ we write
$x\geqq_{P}x^{\prime}$ for $x-x^{\prime}\in P$, $x\geq_{P}x^{\prime}$ for
$x-x^{\prime}\in P\setminus\{0\}$ and $x>_{P}x^{\prime}$ for $x-x^{\prime}%
\in\operatorname*{int}P$.

Recall that the domain $\operatorname*{dom}f$ of the function $f:X\rightarrow
\overline{\mathbb{R}}:=\mathbb{R}\cup\{-\infty,+\infty\}$ is the set $\{x\in
X\mid f(x)<\infty\}$, while the epigraph of $f$ is the set
$\operatorname*{epi}f:=\{(x,\lambda)\in X\times\mathbb{R}\mid f(x)\leq
\lambda\}$; $f$ is proper if $\operatorname*{dom}f\neq\emptyset$ and
$f(x)\neq-\infty$ for every $x\in X$. If $f$ is proper, its recession function
is $f_{\infty}:X\rightarrow\overline{\mathbb{R}}$ whose epigraph is
$(\operatorname*{epi}f)_{\infty}$; $f_{\infty}$ is lower semi\-continuous (lsc
for short) and positively homogeneous, that is $f_{\infty}(tu)=tf_{\infty}(u)$
for all $u\in X$ and $t\in\mathbb{P}:=(0,\infty)$. The conjugate of $f$ is the
function
\[
f^{\ast}:X^{\ast}\rightarrow\overline{\mathbb{R}},\quad f^{\ast}(x^{\ast
}):=\sup\left\{  \left\langle x,x^{\ast}\right\rangle -f(x)\mid x\in
X\right\}
\]
and the subdifferential of the proper function $f$ at $x\in\operatorname*{dom}%
f$ is
\[
\partial f(x):=\left\{  x^{\ast}\in X^{\ast}\mid\left\langle x^{\prime
}-x,x^{\ast}\right\rangle \leq f(x^{\prime})-f(x)\ \forall x^{\prime}\in
X\right\}
\]
and $\partial f(x)=\emptyset$ for $x\in X\setminus\operatorname*{dom}f$. Of
course, the domain of $\partial f$ is $\operatorname*{dom}\partial f:=\{x\in
X\mid\partial f(x)\neq\emptyset\}$ $(\subset\operatorname*{dom}f)$. Clearly,
$\sigma_{A}=(\iota_{A})^{\ast}$, where the indicator function $\iota_{A}$ of
$A\subset X$ is defined by $\iota_{A}(x):=0$ for $x\in A$ and $\iota
_{A}(x):=+\infty$ for $x\in X\setminus A.$

Coming back to the support function, it is well known that for the nonempty
set $A\subset X$ we have that $\sigma_{A}=\sigma_{\operatorname*{conv}%
A}=\sigma_{\operatorname*{cl}A}$, which shows that (in many problems) we can
assume that $A$ is a (nonempty) closed convex set.

For $\emptyset\neq C=\overline{\operatorname*{conv}}\,C$ (that is $C$ is a
nonempty closed convex set) we have that
\[
\left(  C_{\infty}\right)  ^{-}=\operatorname*{cl}(\operatorname*{dom}%
\sigma_{C}),
\]
whence%
\begin{equation}
\operatorname*{int}(\operatorname*{dom}\sigma_{C})=\operatorname*{int}\left(
C_{\infty}\right)  ^{-}. \label{r0}%
\end{equation}
Hence $\operatorname*{int}(\operatorname*{dom}\sigma_{C})\neq\emptyset$ if and
only if $C_{\infty}$ is a pointed cone. Moreover,%
\[
\operatorname*{lin}\nolimits_{0}C=X\Leftrightarrow(\operatorname*{lin}%
\nolimits_{0}C)^{\perp}=\{0\}\Leftrightarrow\operatorname*{int}C\neq
\emptyset,
\]
and $\partial\sigma_{C}(x^{\ast})=C$ for every $x^{\ast}\in
(\operatorname*{lin}_{0}C)^{\perp}$, whence $(\operatorname*{lin}%
\nolimits_{0}C)^{\perp}\subset\operatorname*{dom}\partial\sigma_{C}$. It
follows that $\sigma_{C}$ is differentiable at $x^{\ast}\in
(\operatorname*{lin}_{0}C)^{\perp}$ iff $C$ is a singleton (in which case
$(\operatorname*{lin}_{0}C)^{\perp}=\operatorname*{dom}\sigma_{C}=X^{\ast}$
and $\sigma_{C}$ is differentiable). Furthermore,%
\[
(\operatorname*{lin}\nolimits_{0}C)^{\perp}=\operatorname*{dom}\partial
\sigma_{C}\Leftrightarrow(\operatorname*{lin}\nolimits_{0}C)^{\perp
}=\operatorname*{dom}\sigma_{C}\Leftrightarrow C=\operatorname*{aff}C.
\]
If $C$ is unbounded (or, equivalently, $C$ is not compact) then
$(\operatorname*{lin}_{0}C)^{\perp}\cap\operatorname*{int}(\operatorname*{dom}%
\sigma_{C})=\emptyset.$

In \cite[Prop.\ 1]{Zal:12} it is shown that for the nonempty closed convex set
$C\subset X$ one has that $C_{\infty}$ is pointed iff there are $\overline
{x}\in X$ and a closed convex pointed cone $P$ such that $C\subset\overline
{x}+P$. From this we get immediately that for the nonempty set $A\subset X$
one has that $\operatorname*{int}(\operatorname*{dom}\sigma_{A})\neq\emptyset$
iff there exist $\overline{x}\in X$ and a closed convex pointed cone $P$ such
that $A\subset\overline{x}+P$. Because $\operatorname*{dom}\sigma
_{A+x}=\operatorname*{dom}\sigma_{A}$ for every $x\in X$, we may (and we shall
often do) assume that $A\subset P.$

We have that
\begin{equation}
\partial\sigma_{A}(0)=\overline{\operatorname*{conv}}\,A,\quad\partial
\sigma_{A}(x^{\ast})=\{u\in\overline{\operatorname*{conv}}\,A\mid\left\langle
u,x^{\ast}\right\rangle =\sigma_{A}(x^{\ast})\} \label{r-sub}%
\end{equation}
for every $x^{\ast}\in X^{\ast}$.

Because $\sigma_{A}$ is a sublinear (hence convex) function,
$\sigma_{A}$ is locally Lipschitz on the interior of its domain, and
so its G\^{a}teaux and Fr\'{e}chet differentiability coincide. This is the
reason for speaking simply about the differentiability of
$\sigma_{A}$ in the sequel.

Theorem 25.1 in \cite{Roc:70} states that the proper convex function
$f:\mathbb{R}^{n}\rightarrow\overline{\mathbb{R}}$ is differentiable
at $\overline{x}\in\operatorname*{dom}f$ if and only if $\partial
f(\overline {x})$ is a singleton, in which case $\partial
f(\overline{x})=\left\{  \nabla f(\overline{x})\right\}  $. Hence
$\sigma_{A}$ is differentiable at $\overline{x}^{\ast}$ if and only
if $\partial\sigma_{A}(\overline{x}^{\ast})$ is a singleton. Using
(\ref{r-sub}) it follows that $\sigma_{A}$ is differentiable at $0$
if and only if $A$ is a singleton, in which case $\sigma_{A}$ is a
linear functional.

In the next section we recall some results related to the differentiability of
$\sigma_{A}$ in the case $A$ is a closed convex set with $\operatorname*{int}%
(\operatorname*{dom}\sigma_{A})\neq\emptyset$. These results suggest the kind
of conditions to be imposed in order that the differentiability of $\sigma
_{A}$ imply the convexity of $A.$

\section{The convex case}

Let $C\subset X$ be a nonempty closed convex set with $\operatorname*{int}%
(\operatorname*{dom}\sigma_{C})\neq\emptyset$, or equivalently, there exist
$\overline{x}\in X$ and a pointed closed convex cone $P\subset X$ with
$C\subset\overline{x}+P$. We recall some results concerning the
differentiability of $\sigma_{C}$ which can be found in \cite{Zal:12}.

\begin{theorem}
\label{ps1}\emph{(\cite[Thm.\ 1]{Zal:12})} $\sigma_{C}$ is
differentiable on
$\operatorname*{dom}\partial\sigma_{C}\setminus(\operatorname*{lin}%
_{0}C)^{\perp}$ if and only if
\begin{equation}
\forall x,x^{\prime}\in C,\ x\neq x^{\prime},\ \forall\lambda\in(0,1):\lambda
x+(1-\lambda)x^{\prime}\in\operatorname*{rint}C, \label{rs1a}%
\end{equation}
or, equivalently,%
\begin{equation}
\forall x,x^{\prime}\in\operatorname*{rbd}C,\ x\neq x^{\prime},\ \forall
\lambda\in(0,1):\lambda x+(1-\lambda)x^{\prime}\notin\operatorname*{rbd}C.
\label{r-s1a}%
\end{equation}

\end{theorem}

Recall that a closed convex set $C$ with nonempty interior (hence
$\operatorname*{rint}C=\operatorname*{int}C$) is called strictly convex if
condition (\ref{rs1a}) is verified.

In the case in which $C$ is compact the following result holds.

\begin{corollary}
\label{c-compact}\emph{(\cite[Cor.\ 4]{Zal:12})} Assume that $C$ is compact. Then

\emph{(i)} $\sigma_{C}$ is differentiable on $X^{\ast}\setminus
(\operatorname*{lin}_{0}C)^{\perp}$ if and only if $C$ verifies condition
$(\ref{rs1a})$.

\emph{(ii)} $\sigma_{C}$ is differentiable on $X^{\ast}\setminus\{0\}$ if and
only if either $C$ is a singleton or $\operatorname*{int}C\neq\emptyset$ and
$C$ is strictly convex.
\end{corollary}

In the case in which $C$ is not compact one has the next result.

\begin{theorem}
\label{ps3}\emph{(\cite[Thm.\ 6]{Zal:12})} Assume that $C$ is
unbounded. Then
$\sigma_{C}$ is differentiable on $\operatorname*{int}(\operatorname*{dom}%
\sigma_{C})$ if and only if
\begin{equation}
\forall x,x^{\prime}\in S_{C},\ x\neq x^{\prime},\ \forall\lambda
\in(0,1):\lambda x+(1-\lambda)x^{\prime}\notin S_{C}, \label{fps-z}%
\end{equation}
where
\begin{equation*}
S_{C}:=\partial\sigma_{C}(\operatorname*{int}(\operatorname*{dom}\sigma_{C})).
\label{r-sc}%
\end{equation*}

\end{theorem}

One has also the following result.

\begin{proposition}
\label{ps2}\emph{(\cite[Prop.\ 7]{Zal:12})} Assume that $C$ is unbounded. If
\begin{equation}
\forall x,x^{\prime}\in E_{C},\ x\neq x^{\prime},\ \forall\lambda
\in(0,1):\lambda x+(1-\lambda)x^{\prime}\notin E_{C}, \label{rs2c}%
\end{equation}
where
\begin{equation*}
E_{C}:=C\setminus\left[  C+(C_{\infty}\setminus\{0\})\right]  , 
\end{equation*}
or, equivalently,%
\begin{equation}
\forall x,x^{\prime}\in C,\ x\neq x^{\prime},\ \forall\lambda\in(0,1):\lambda
x+(1-\lambda)x^{\prime}\in C+(C_{\infty}\setminus\{0\}), \label{rs2}%
\end{equation}
then $\sigma_{C}$ is differentiable on $\operatorname*{int}%
(\operatorname*{dom}\sigma_{C})$. Moreover, if $\dim(\operatorname*{lin}%
_{0}C)\leq2$ then the converse is also true.
\end{proposition}

\begin{remark}\emph{
\label{rem1}In \cite[Lem.\ 5]{Zal:12} it is shown that for $C$ unbounded one
has $S_{C}\subset E_{C}\subset\operatorname*{rbd}C$, and so $(\ref{r-s1a}%
)$~$\Rightarrow$ $(\ref{rs2c})$~$\Leftrightarrow$~$(\ref{rs2})$~$\Rightarrow$
$(\ref{fps-z})$. In fact, taking into account that $\operatorname*{int}%
(\operatorname*{dom}\sigma_{C})\subset\operatorname*{dom}\partial\sigma_{C},$
if $C$ is unbounded and $\sigma_{C}$ is differentiable on $\operatorname*{dom}%
\partial\sigma_{C}\setminus(\operatorname*{lin}_{0}C)^{\perp}$, then
$\operatorname*{dom}\partial\sigma_{C}\setminus(\operatorname*{lin}%
_{0}C)^{\perp}=\operatorname*{int}(\operatorname*{dom}\sigma_{C}).$}
\end{remark}

In \cite{Zal:12} it is shown that the set $A$ defined in (\ref{d4}) is a
closed convex set included in $\mathbb{R}_{+}^{3}$ with $A_{\infty}%
=\mathbb{R}_{+}^{3}$ for which $\sigma_{A}$ is differentiable on
$\operatorname*{int}(\operatorname*{dom}\sigma_{A})$ but for which (\ref{rs2})
does not hold (\cite[Prop.\ 10]{Zal:12}). Hence, in general, (\ref{rs2}) does
not imply (\ref{r-s1a}).

\begin{theorem}
\label{Fact 7}\emph{(\cite[Thm.\ 12]{Zal:12})} Let $K\subset X$ be a
pointed closed convex cone and let $A\subset X$ be a nonempty closed
convex set such that $A_{\infty}\subset K$. Then $\sigma_{A}$ is
differentiable on $-K^{\#}$
if and only if%
\begin{equation}
\forall x,x^{\prime}\in SE(A;K),\text{ }x\neq x^{\prime},\ \forall\lambda
\in(0,1):\lambda x+(1-\lambda)x^{\prime}\notin SE(A;K), \label{fps}%
\end{equation}
where%
\begin{equation*}
SE(A;K):=\cup_{x^{\ast}\in K^{\#}}\partial\sigma_{A}(-x^{\ast})=\partial
\sigma_{A}(-K^{\#}). 
\end{equation*}

\end{theorem}

Note that Theorem \ref{ps3} can be obtained from Theorem \ref{Fact 7} taking
$K:=A_{\infty}$ because $S_{A}=SE(A;A_{\infty})$ and $K^{\#}%
=-\operatorname*{int}(\operatorname*{dom}\sigma_{A}).$

\section{Relations between the differentiability of $\sigma_{A}$ and the
convexity of $A$}

Throughout this section $\mathbb{N}\ni p:= \dim X\geq2$ (the case $\dim X=1$
being trivial) and $A\subset X$ is a nonempty closed set with the property
that $P:=(\overline{\operatorname*{conv}}\,A)_{\infty}$ is pointed. We
consider the multi\-function
\begin{equation}
W_{A}:X^{\ast}\rightrightarrows X,\quad W_{A}(x^{\ast}):=\{a\in A\mid
\left\langle a,x^{\ast}\right\rangle =\sigma_{A}(x^{\ast})\}=A\cap
\partial\sigma_{A}(x^{\ast}). \label{r1}%
\end{equation}
Of course, $\operatorname*{dom}W_{A}\subset\operatorname*{dom}\sigma
_{A}\subset P^{-}$.

\begin{proposition}
\label{Fact 1}The following assertions hold:

\emph{(a)}~$\operatorname*{int}P^{-}=\operatorname*{int}(\operatorname*{dom}%
\sigma_{A})\neq\emptyset$ and $\sigma_{A}$ is continuous on
$\operatorname*{int}P^{-}$;

\emph{(b)}~$W_{A}(x^{\ast})$ is nonempty and compact for every $x^{\ast}%
\in\operatorname*{int}P^{-}$;

\emph{(c)}~$\partial\sigma_{A}(x^{\ast})$ is nonempty convex and compact for
every $x^{\ast}\in\operatorname*{int}P^{-}.$
\end{proposition}

Proof. (a), (c) The equality $\operatorname*{int}P^{-}=\operatorname*{int}%
(\operatorname*{dom}\sigma_{A})$ was observed above. Because $\sigma_{A}$ is
convex we have that $\sigma_{A}$ is continuous on $\operatorname*{int}%
(\operatorname*{dom}\sigma_{A})$ by a well-known result in Convex Analysis.
Since the subdifferential of a proper convex function is nonempty, convex and
compact at any point of continuity from its domain, (c) follows.

(b) Fix $x^{\ast}\in\operatorname*{int}P^{-}=-P^{\#}$. Let $u_{0}\in A$ be
fixed and take $A_{0}:=\{u\in A\mid\left\langle u,x^{\ast}\right\rangle
\geq\left\langle u_{0},x^{\ast}\right\rangle \}$. Clearly $A_{0}$ is nonempty
and closed. Assume that $A_{0}$ is not bounded. Then there exists
$(u_{n})\subset A_{0}$ with $\left\Vert u_{n}\right\Vert \rightarrow\infty$.
We may assume that $\left\Vert u_{n}\right\Vert ^{-1}u_{n}\rightarrow v$, and
so $v\in P\setminus\{0\}$. Since $\left\langle u_{n},x^{\ast}\right\rangle
\geq\left\langle u_{0},x^{\ast}\right\rangle $ for every $n$, dividing by
$\left\Vert u_{n}\right\Vert $ and passing to the limit we get the
contradiction $0>\left\langle v,x^{\ast}\right\rangle \geq0$. Hence $A_{0}$ is
bounded, and so $A_{0}$ is compact. Since $\sigma_{A}(x^{\ast})=\sup\left\{
\left\langle u,x^{\ast}\right\rangle \mid x^{\ast}\in A_{0}\right\}  $ and
$A_{0}$ is compact, there exists $u\in A_{0}\subset A$ with $\left\langle
u,x^{\ast}\right\rangle =\sigma_{A}(x^{\ast})\in\mathbb{R}$. It follows that
$W_{A}(x^{\ast})$ is nonempty and compact. \hfill$\square$

\medskip

Note that we can have that $W_{A}(x^{\ast})$ is nonempty and compact without
having $x^{\ast}\in\operatorname*{int}P^{-}.$

\begin{example}\emph{
\label{ex1}Take $X=\mathbb{R}^{2}$ endowed with the Euclidean norm and
\[
A:=\left\{  (a,b)\in\mathbb{R}^{2}\mid b\geq\left\vert a\right\vert \left(
1+(a^{2}+1)^{-1}\right)  \right\}  .
\]
Then $C:=\overline{\operatorname*{conv}}\,A=\left\{  (a,b)\in\mathbb{R}%
^{2}\mid b\geq\left\vert a\right\vert \right\}  $. We have that
$\operatorname*{dom}\sigma_{A}=-C$ and $W_{A}(-1,-1)=\{(0,0)\}$ is compact and
nonempty; clearly, $(0,0)\notin\operatorname*{int}(\operatorname*{dom}%
\sigma_{A}).$}
\end{example}

However, the next result holds.

\begin{proposition}
\label{Fact 5}Let $B\subset X$ be a nonempty closed set. Then $x^{\ast}%
\in\operatorname*{int}(\operatorname*{dom}\sigma_{B})$ if and only if
$\partial\sigma_{B}(x^{\ast})$ is nonempty and compact.
\end{proposition}

Proof. Set $C:=\overline{\operatorname*{conv}}\,B\ $and $Q:=C_{\infty}$. If
$x^{\ast}\in\operatorname*{int}(\operatorname*{dom}\sigma_{B})$ then, by
(\ref{r0}), $\operatorname*{int}Q^{-}\neq\emptyset$, and so $Q^{-}$ is
pointed. From Proposition \ref{Fact 1}~(c) we have that $\partial\sigma
_{B}(x^{\ast})$ is nonempty and compact.

Assume now that $\partial\sigma_{B}(x^{\ast})$ is nonempty and compact. Take
$f:X\rightarrow\overline{\mathbb{R}}$, $f:=-x^{\ast}+\iota_{C}$; then $f$ is a
proper lsc convex function, and $\{x\in X\mid f(x)=\inf f\}=\partial\sigma
_{B}(x^{\ast})$. Since $\sigma_{B}(x^{\ast})$ is nonempty and compact, it
follows that $0\in\operatorname*{int}(\operatorname*{dom}f^{\ast})$, that is,
$x^{\ast}\in\operatorname*{int}(\operatorname*{dom}\sigma_{C})$ (see e.g.
\cite[Exer.\ 2.41]{Zal:02}). \hfill$\square$

\begin{proposition}
\label{Fact 3}One has that $\partial\sigma_{A}(x^{\ast})=\operatorname*{conv}%
W_{A}(x^{\ast})\subset\operatorname*{conv}A$ for every $x^{\ast}%
\in\operatorname*{int}P^{-}$.
\end{proposition}

Proof. Fix $x^{\ast}\in\operatorname*{int}P^{-}$. From (\ref{r1}) we have that
$W_{A}(x^{\ast})\subset\partial\sigma_{A}(x^{\ast})$, and so
$\operatorname*{conv}W_{A}(x^{\ast})\subset\partial\sigma_{A}(x^{\ast})$. Let
$u\in\partial\sigma_{A}(x^{\ast})$. Hence $u\in\overline{\operatorname*{conv}%
}\,A$ and $\left\langle x^{\ast},u\right\rangle =\sigma_{A}(x^{\ast})$. Using
the Carath\'{e}odory theorem, we find $(\lambda_{n}^{k})_{n\geq1}%
\subset\lbrack0,1]$ and $(u_{n}^{k})_{n\geq1}\subset A$ for $k\in
\overline{1,p+1}$ such that $\sum_{k=1}^{p+1}\lambda_{n}^{k}=1$, $u_{n}%
:=\sum_{k=1}^{p+1}\lambda_{n}^{k}u_{n}^{k}\rightarrow u$.

We claim that the sequences $(\lambda_{n}^{k}u_{n}^{k})_{n\geq1}$ are bounded.
In the contrary case we may assume that
\[
\left\Vert \lambda_{n}^{1}u_{n}^{1}\right\Vert \leq\left\Vert \lambda_{n}%
^{2}u_{n}^{2}\right\Vert \leq\cdots\leq\left\Vert \lambda_{n}^{p+1}u_{n}%
^{p+1}\right\Vert \quad\forall n\geq1
\]
and $\big\Vert \lambda_{n}^{p+1}u_{n}^{p+1}\big\Vert
\rightarrow\infty$ as $n\rightarrow\infty$. Taking sub\-sequences if
necessary, we may assume that $\big\Vert \lambda_{n}^{p+1}u_{n}^{p+1}%
\big\Vert ^{-1}\lambda_{n}^{k}u_{n}^{k}\rightarrow v_{k}$ for $k\in
\overline{1,p+1}$. Since $\lambda_{n}^{k}\in\lbrack0,1]$ and
$\big\Vert \lambda_{n}^{p+1}u_{n}^{p+1}\big\Vert
\rightarrow\infty$ we have that $\big\Vert \lambda_{n}^{p+1}u_{n}%
^{p+1}\big\Vert ^{-1}\lambda_{n}^{k}\rightarrow0$. Since $u_{n}^{k}\in
A\subset\overline{\operatorname*{conv}}\,A$, we obtain that $v_{k}%
\in(\overline{\operatorname*{conv}}\,A)_{\infty}=P$ for every $k\in
\overline{1,p+1}$. From $\sum_{k=1}^{p+1}\lambda_{n}^{k}u_{n}^{k}\rightarrow
u$ we get $v_{1}+\cdots+v_{p+1}=0$. Since $P$ is pointed and $v_{p+1}\neq0$ we
get a contradiction. Hence the sequences $(\lambda_{n}^{k}u_{n}^{k})_{n\geq1}$
are bounded is true (for $k\in\overline{1,p+1}$).

We may assume that $\lambda_{n}^{k}\rightarrow\lambda^{k}$ and $\lambda
_{n}^{k}u_{n}^{k}\rightarrow v^{k}$ for every $k\in\overline{1,p+1}$;
moreover, we may assume that $\lambda^{1}\leq\lambda^{2}\leq\ldots\leq
\lambda^{p+1}$. Since $\sum_{k=1}^{p+1}\lambda^{k}=1$, $\lambda^{p+1}>0$. If
$\lambda^{k}>0$ then $u_{n}^{k}\rightarrow u^{k}:=(\lambda^{k})^{-1}v_{k}\in
A$. If $\lambda^{1}>0$ then $u=\sum_{k=1}^{p+1}\lambda^{k}u^{k}\in
\operatorname*{conv}A$. If $\lambda^{1}=0$ take $k_{0}\in\overline{1,p}$ such
that $\lambda^{k_{0}}=0$ and $\lambda^{k_{0}+1}>0$. Then $u=\overline
{v}+\overline{u}$ with $\overline{v}:=\sum_{k=1}^{k_{0}}v^{k}$, $\overline
{u}:=\sum_{k=k_{0}+1}^{p+1}\lambda^{k}u^{k}\in\operatorname*{conv}A$. Since
$\lambda^{k}=0$ for $1\leq k\leq k_{0}$, we have that $v_{k}\in A_{\infty
}\subset P$ for such $k$, and so $\overline{v}\in P$. We have that
\[
\sigma_{A}(x^{\ast})=\left\langle u,x^{\ast}\right\rangle =\left\langle
\overline{v},x^{\ast}\right\rangle +\left\langle \overline{u},x^{\ast
}\right\rangle \leq\left\langle \overline{u},x^{\ast}\right\rangle
=\sum_{k=k_{0}+1}^{p+1}\lambda^{k}\big\langle u^{k},x^{\ast}\big\rangle\leq
\sigma_{A}(x^{\ast}).
\]
Therefore, $\left\langle \overline{v},x^{\ast}\right\rangle =0$, and so
$\overline{v}=0$ because $x^{\ast}\in\operatorname*{int}P^{+}=P^{\#}$. It
follows that $u=\overline{u}\in\operatorname*{conv}A$. \hfill$\square$

\begin{proposition}
\label{Fact 4}Let $x^{\ast}\in\operatorname*{int}P^{-}$. Then $\sigma_{A}$ is
differentiable at $x^{\ast}$ if and only if $W_{A}(x^{\ast})$ is a singleton.
In this case $\nabla\sigma_{A}(x^{\ast})\in A$ and $W_{A}(x^{\ast}%
)=\{\nabla\sigma_{A}(x^{\ast})\}.$
\end{proposition}

Proof. Because $\sigma_{A}$ is convex and continuous at $x^{\ast}%
\in\operatorname*{int}(\operatorname*{dom}\sigma_{A})$, $\sigma_{A}$
is differentiable at $x^{\ast}$ if and only if $\partial\sigma
_{A}(x^{\ast})$ is a singleton. Using Proposition \ref{Fact 3}, this
happens exactly when $W_{A}(x^{\ast})$ is a singleton.
\hfill$\square$

\begin{corollary}
\label{cor2}If $W_{A}(x^{\ast})$ is a singleton for every $x^{\ast}%
\in\operatorname*{int}P^{-}$ then $\sigma_{A}$ is differentiable on
$\operatorname*{int}(\operatorname*{dom}\sigma_{A})=\operatorname*{int}P^{-}.$
\end{corollary}

\begin{corollary}
\label{cor3}If $\sigma_{A}$ is differentiable on $\operatorname*{dom}%
\partial\sigma_{A}\setminus(\operatorname*{lin}_{0}A)^{\perp}$ then
$\operatorname*{rbd}(\overline{\operatorname*{conv}}\,A)\subset A.$
\end{corollary}

Proof. Set $C:=\overline{\operatorname*{conv}}\,A$; it follows that
$\operatorname*{lin}_{0}A=\operatorname*{lin}_{0}C$.

Let $u\in\operatorname*{rbd}C$. Using a separation theorem, there exists
$x^{\ast}\in X^{\ast}\setminus(\operatorname*{lin}_{0}C)^{\perp}$ (that is
$x^{\ast}$ is not constant on $C$) such that $\left\langle u,x^{\ast
}\right\rangle =\sigma_{C}(x^{\ast})=\sigma_{A}(x^{\ast})$. Therefore,
$x^{\ast}\in\operatorname*{dom}\partial\sigma_{A}\setminus(\operatorname*{lin}%
_{0}A)^{\perp}$ and $u\in\partial\sigma_{A}(x^{\ast})$. Because $\sigma_{A}$
is differentiable at $x^{\ast}$, by Proposition \ref{Fact 4} we obtain that
$\nabla\sigma_{A}(x^{\ast})=u\in A$. \hfill$\square$

\medskip

Note that we cannot obtain the convexity of $A$ in Corollary \ref{cor3} under
its hypothesis.

\begin{example}\label{ex2}\emph{
Consider $X:=\mathbb{R}^{2}$ endowed with the Euclidean norm and the
sets $A_{1}:=\{(x,y)\in X\mid\left\Vert (x,y)\right\Vert =1\}$,
$A_{2}:=\{x\in X\mid1/2\leq\left\Vert (x,y)\right\Vert \leq1\},$
$A_{3}:=\{(x,y)\in X\mid x>0,\ y=1/x\}$ and $A_{4}:=\{(x,y)\in X\mid
x>0,\ 1/x\leq y\leq2/x\}$. Then
$\sigma_{A_{1}}(u,v)=\sigma_{A_{2}}(u,v)=\left\Vert (u,v)\right\Vert
=\sqrt{u^{2}+v^{2}}$ and
$\sigma_{A_{3}}(u,v)=\sigma_{A_{4}}(u,v)=-2\sqrt{uv}$ for
$u,v\leq0$, $\sigma_{A_{3}}(u,v)=\sigma_{A_{4}}(u,v)=+\infty$
otherwise.
Clearly, $\sigma_{A_{i}}$ is differentiable on $\operatorname*{dom}%
\partial\sigma_{A_{i}}\setminus(\operatorname*{lin}_{0}A_{i})^{\perp
}=\operatorname*{int}(\operatorname*{dom}\sigma_{A_{i}})\setminus\{0\}$.
Note that $A_{i}=\operatorname*{cl}(\operatorname*{int}A_{i})$ for
$i\in\{2,4\}.$}
\end{example}

These simple examples show that there is no hope to get the convexity of $A$
from the differentiability of $\sigma_{A}$ on $\operatorname*{dom}%
\partial\sigma_{A}\setminus(\operatorname*{lin}_{0}A)^{\perp}$ or on
$\operatorname*{int}(\operatorname*{dom}\sigma_{A})\setminus\{0\}$ in the case
in which $A$ is bounded. Even for $A$ unbounded one needs supplementary
conditions. In the sequel we concentrate on the case in which $A$ is unbounded.

In Theorem \ref{Fact 14} below we provide a supplementary condition on $A$ to
be added in Corollary \ref{cor3} in order to get the convexity of $A$. First
we establish an auxiliary result.

\begin{lemma}
\label{Fact 13}Assume that $A_{\infty}$ is a pointed convex cone and
$A=A+A_{\infty}$. Then $\operatorname*{conv}A$ is closed.
\end{lemma}

Proof. Let $u\in\overline{\operatorname*{conv}}\,A$. Using the
Carath\'{e}odory theorem, we find $(\lambda_{n}^{k})_{n\geq1}\subset
\lbrack0,1]$ and $(u_{n}^{k})_{n\geq1}\subset A$ for $k\in\overline{1,p+1}$
such that $\sum_{k=1}^{p+1}\lambda_{n}^{k}=1$, $u_{n}:=\sum_{k=1}^{p+1}%
\lambda_{n}^{k}u_{n}^{k}\rightarrow u$. As in the proof of Proposition
\ref{Fact 3}, we have that the sequences $(\lambda_{n}^{k}u_{n}^{k})_{n\geq1}$
are bounded. In the contrary case we may assume that
\[
\left\Vert \lambda_{n}^{1}u_{n}^{1}\right\Vert \leq\left\Vert \lambda_{n}%
^{2}u_{n}^{2}\right\Vert \leq\cdots\leq\left\Vert \lambda_{n}^{p+1}u_{n}%
^{p+1}\right\Vert \quad\forall n\geq1
\]
and $\big\Vert\lambda_{n}^{p+1}u_{n}^{p+1}\big\Vert\rightarrow\infty$ as
$n\rightarrow\infty$. Taking sub\-sequences if necessary, we may assume that
$\big\Vert\lambda_{n}^{p+1}u_{n}^{p+1}\big\Vert^{-1}\lambda_{n}^{k}u_{n}%
^{k}\rightarrow v_{k}$ for $k\in\overline{1,p+1}$. Since $\lambda_{n}^{k}%
\in\lbrack0,1]$ and $\big\Vert\lambda_{n}^{p+1}u_{n}^{p+1}\big\Vert\rightarrow
\infty$ we have that $\big\Vert\lambda_{n}^{p+1}u_{n}^{p+1}\big\Vert^{-1}%
\lambda_{n}^{k}\rightarrow0$. Since $u_{n}^{k}\in A$, we obtain that $v_{k}\in
A_{\infty}=:K$ for every $k\in\overline{1,p+1}$. From $\sum_{k=1}^{p+1}%
\lambda_{n}^{k}u_{n}^{k}\rightarrow u$ we get $v_{1}+\cdots+v_{p+1}=0$. Since
$K$ is convex we get $v_{1}+\cdots+v_{p}=-v_{p+1}\in K\cap(-K)$. Because $K$
is pointed we get the contradiction $v_{p+1}=0$. Hence the sequences
$(\lambda_{n}^{k}u_{n}^{k})_{n\geq1}$ are bounded. We may assume that
$\lambda_{n}^{k}\rightarrow\lambda^{k}$ and $\lambda_{n}^{k}u_{n}%
^{k}\rightarrow v^{k}\in X$ for every $k\in\overline{1,p+1}$; moreover, we may
assume that $\lambda^{1}\leq\lambda^{2}\leq\ldots\leq\lambda^{p+1}$. Since
$\sum_{k=1}^{p+1}\lambda^{k}=1$, $\lambda^{p+1}>0$. If $\lambda^{k}>0$ then
$u_{n}^{k}\rightarrow u^{k}:=(\lambda^{k})^{-1}v_{k}\in A$. If $\lambda^{1}>0$
then $u=\sum_{k=1}^{p+1}\lambda^{k}u^{k}\in\operatorname*{conv}A$. If
$\lambda^{1}=0$ take $k_{0}\in\overline{1,p}$ such that $\lambda^{k_{0}}=0$
and $\lambda^{k_{0}+1}>0$. Then $u=\overline{v}+\overline{u}$ with
$\overline{v}:=\sum_{k=1}^{k_{0}}v^{k}\in K$, $\overline{u}:=\sum_{k=k_{0}%
+1}^{p+1}\lambda^{k}u^{k}$, and so $u=\sum_{k=k_{0}+1}^{p+1}\lambda^{k}%
(u^{k}+\overline{v})\in\operatorname*{conv}A$. It follows that $\overline
{\operatorname*{conv}}\,A\subset\operatorname*{conv}A$. \hfill$\square$

\medskip

A nice application of the preceding lemma is the fact that the convex hull of
the epigraph of a proper lsc $1$-coercive function $f:\mathbb{R}%
^{p}\rightarrow\overline{\mathbb{R}}$ is closed, result which can be found in
\cite{Val:70}, \cite{HirLem:93a}, \cite{BenHir:96}. First we give the next
result which is probably known.

\begin{proposition}
\label{prop-rec}Let $f:X\rightarrow\overline{\mathbb{R}}$ be a proper lsc
function. Then $f$ is $1$-coercive (that is $\lim_{\left\Vert x\right\Vert
\rightarrow\infty}f(x)/\left\Vert x\right\Vert =\infty$) if and only if
$f_{\infty}=\iota_{\{0\}}$.
\end{proposition}

Proof. Assume first that $\lim_{\left\Vert x\right\Vert \rightarrow\infty
}f(x)/\left\Vert x\right\Vert =\infty$ and take $(u,\alpha)\in
(\operatorname*{epi}f)_{\infty}$. Then there exist the sequences
$(t_{n})_{n\geq1}\subset\mathbb{P}$ and $\left(  (x_{n},\lambda_{n})\right)
_{n\geq1}\subset\operatorname*{epi}f$ such that $t_{n}\rightarrow0$ and
$t_{n}(x_{n},\lambda_{n})\rightarrow(u,\alpha)$. Suppose that $u\neq0$. Then
$\left\Vert x_{n}\right\Vert \rightarrow\infty$, and so $f(x_{n})/\left\Vert
x_{n}\right\Vert \rightarrow\infty$. Since $f(x_{n})/\left\Vert x_{n}%
\right\Vert \leq(t_{n}\lambda_{n})/\left\Vert t_{n}x_{n}\right\Vert
\rightarrow\alpha/\left\Vert u\right\Vert $, we get the contradiction
$\infty\leq\alpha/\left\Vert u\right\Vert $. Hence $u=0$. If $\alpha<0$ then
$t_{n}\lambda_{n}\leq\alpha/2<0$ for large $n$, and so $f(x_{n})\leq
\alpha/(2t_{n})$ for such $n$. Hence $f(x_{n})\rightarrow-\infty$. Because $f$
is $1$-coercive, we obtain that $(x_{n})$ is bounded, and so, passing to a
subsequence if necessary, we assume that $x_{n}\rightarrow x$. Since $f$ is
lsc, we get the contradiction $-\infty<f(x)\leq\liminf f(x_{n})=-\infty.$

Assume now that $\liminf_{\left\Vert x\right\Vert \rightarrow\infty
}f(x)/\left\Vert x\right\Vert <\infty$. Then there exist a sequence
$(x_{n})_{n\geq1}\subset\mathbb{R}_{+}$ and $\alpha\in\mathbb{R}$
such that $\left\Vert x_{n}\right\Vert \rightarrow\infty$ and
$f(x_{n})/\left\Vert x_{n}\right\Vert \leq\alpha$ for every $n$.
Passing to a subsequence if necessary, we assume that
$x_{n}/\left\Vert x_{n}\right\Vert \rightarrow u$ $(\neq0)$. Since
$(x_{n},\alpha\left\Vert x_{n}\right\Vert )\in
\operatorname*{epi}f$ and $0<t_{n}:=\left\Vert x_{n}\right\Vert ^{-1}%
\rightarrow0$, we get $(u,\alpha)=\lim t_{n}(x_{n},\alpha\left\Vert
x_{n}\right\Vert )\in(\operatorname*{epi}f)_{\infty}$. Hence
$(\operatorname*{epi}f)_{\infty}\not =\{0\}\times\mathbb{R}_{+}$.
The proof is complete. \hfill$\square$

\begin{corollary}
\label{cor-bhu}Let $f:X\rightarrow\overline{\mathbb{R}}$ be a proper lsc
function. Assume that $f_{\infty}$ is a proper convex (hence sublinear)
function such that $f_{\infty}(u)+f_{\infty}(-u)>0$ for all $u\in
X\setminus\{0\}$ and $f(x+u)\leq f(x)+f_{\infty}(u)$ for all $x,u\in X$; then
$\operatorname*{conv}(\operatorname*{epi}f)$ is closed. In particular,
$\operatorname*{conv}(\operatorname*{epi}f)$ is closed if $f$ is $1$-coercive.
\end{corollary}

Proof. Since $f_{\infty}$ is convex we have that $(\operatorname*{epi}%
f)_{\infty}$ is a closed convex cone. Take $(u,\alpha)\in(\operatorname*{epi}%
f)_{\infty}\cap-(\operatorname*{epi}f)_{\infty}$. Then
$0=f_{\infty}(0)\leq
f_{\infty}(u)+f_{\infty}(-u)\leq\alpha+(-\alpha)=0$, and so $u=0$.
Then $0=f_{\infty}(0)\leq\min\left\{  \alpha,-\alpha\right\}  ,$
whence $\alpha=0$. Hence $(\operatorname*{epi}f)_{\infty}$ is
pointed. Using Lemma \ref{Fact 3} we obtain that
$\operatorname*{conv}(\operatorname*{epi}f)$ is closed.
\hfill$\square$

\medskip

Another example (besides the epigraph of a proper lsc $1$-coercive function)
of set $A$ verifying the hypothesis of Lemma \ref{Fact 13} is when
$A=A+K\subset K$ with $K\subset X$ a proper closed convex pointed cone because
in this case $A_{\infty}=(\overline{\operatorname*{conv}}\,A)_{\infty}=K.$

In the rest of this section we assume that the subsets $A$ and $K$ of $X$
verify the following condition

\begin{itemize}\item[(H)] \textit{$K$ is a pointed closed convex cone with nonempty
interior and $A$ is a closed nonempty set such that $A=A+K\subset
K\setminus\{0\}.$}
\end{itemize}

\noindent In this case $\operatorname*{dom}\sigma_{A}=K^{-}$ and
$\operatorname*{int}(\operatorname*{dom}\sigma_{A})=\operatorname*{int}%
K^{-}=-K^{\#}$. Moreover,
\begin{equation}
\emptyset\neq\operatorname*{int}A=A+\operatorname*{int}K=\operatorname*{int}%
K\cap\operatorname*{int}\nolimits_{K}A\quad\text{and}\quad
A=\operatorname*{cl}(A\cap\operatorname*{int}K); \label{r-intA}%
\end{equation}
it follows that $\operatorname*{int}\nolimits_{K}A=\operatorname*{int}A$, and
so $\operatorname*{bd}_{K}A=\operatorname*{bd}A$, if $A\subset
\operatorname*{int}K.$

\begin{theorem}
\label{Fact 14}Assume that
\begin{equation}
A+(K\setminus\{0\})\subset\operatorname*{int}\nolimits_{K}A. \label{fp-ssc}%
\end{equation}
If $\sigma_{A}$ is differentiable on $\operatorname*{int}K^{-}$ then $A$ is
convex and%
\begin{equation}
\forall a,a^{\prime}\in A\cap\operatorname*{int}K,\ a\neq a^{\prime}%
,\ \forall\lambda\in(0,1):\lambda a+(1-\lambda)a^{\prime}\in
\operatorname*{int}A. \label{r-sa}%
\end{equation}

Conversely, if $\dim X=2$ and $(\ref{r-sa})$ holds, then $A$ is convex and
$(\ref{rs2})$ is verified; therefore, $\sigma_{A}$ is differentiable on
$\operatorname*{int}K^{-}.$
\end{theorem}

Proof. Because $A+(K\setminus\{0\})\subset\operatorname*{int}_{K}A\subset A$,
we have that $A+K=A$. Then $\operatorname*{conv}A=(\operatorname*{conv}A)+K,$
and so $\operatorname*{int}(\operatorname*{conv}A)=(\operatorname*{conv}%
A)+\operatorname*{int}K.$

Assume now that $\sigma_{A}$ is differentiable on $\operatorname*{int}K^{-}$.
As seen above, $A_{\infty}=K$ and $A=A+K$; hence $\operatorname*{conv}A$ is
closed by Lemma \ref{Fact 13}.

Let us prove that $\operatorname*{int}K\cap\operatorname*{conv}A\subset A.$
Take first $\overline{u}\in\operatorname*{int}K\cap\operatorname*{bd}%
(\operatorname*{conv}A)$ and consider $k\in K\setminus\{0\}$. Since
$\overline{u}\in\operatorname*{conv}A$, $\overline{u}=\sum_{i\in I}\lambda
_{i}u_{i}$ for some nonempty finite set $I$, $(\lambda_{i})_{i\in I}%
\subset(0,1)$ with $\sum_{i\in I}\lambda_{i}=1$, and $(u_{i})_{i\in I}\subset
A$. Since $u_{i}+k\in\operatorname*{int}_{K}A$ by our hypothesis, there exists
$\mu_{i}\in(0,1)$ such that $\mu(u_{i}+k)\in A$ for every $\mu\in\lbrack
\mu_{i},1]$. Taking $\overline{\mu}:=\max\{\mu_{i}\mid i\in I\}\in(0,1)$, we
have that $\overline{\mu}(u_{i}+k)\in A$ for every $i\in I$. It follows that
$\overline{\mu}(\overline{u}+k)=\sum_{i\in I}\lambda_{i}\overline{\mu}%
(u_{i}+k)\in\operatorname*{conv}A$, and so
\[
\overline{u}+k=\overline{\mu}(\overline{u}+k)+(1-\overline{\mu})(\overline
{u}+k)\in\operatorname*{conv}A+\operatorname*{int}K\subset\operatorname*{int}%
\left(  \operatorname*{conv}A\right)  .
\]
Hence $\overline{u}+(K\setminus\{0\})\subset\operatorname*{int}%
(\operatorname*{conv}A)$. Since $\overline{u}\in\operatorname*{bd}%
(\operatorname*{conv}A)$, there exists $x^{\ast}\in X^{\ast}\setminus\{0\}$
such that $\sigma_{A}(x^{\ast})=\left\langle \overline{u},x^{\ast
}\right\rangle $. Because $x^{\ast}\neq0$, we have that $\left\langle
\overline{u},x^{\ast}\right\rangle >\left\langle u,x^{\ast}\right\rangle $ for
every $u\in\operatorname*{int}(\operatorname*{conv}A)$. Because $\overline
{u}+(K\setminus\{0\})\subset\operatorname*{int}(\operatorname*{conv}A)$ we get
$\left\langle k,x^{\ast}\right\rangle <0$ for every $k\in K\setminus\{0\}$,
and so $x^{\ast}\in-K^{\#}=\operatorname*{int}(\operatorname*{dom}\sigma_{A}%
)$. Because $\sigma_{A}$ is differentiable at $x^{\ast}$, by Proposition
\ref{Fact 4} we obtain that $\nabla\sigma_{A}(x^{\ast})=\overline{u}\in A.$

Take now $\overline{u}\in\operatorname*{int}K\cap\operatorname*{conv}A$ and
consider $\alpha:=\min\{\gamma>0\mid\gamma\overline{u}\in\operatorname*{conv}%
A\}\in(0,1]$; clearly $\alpha\overline{u}\in\operatorname*{int}K\cap
\operatorname*{bd}(\operatorname*{conv}A)$. By the argument above we have that
$\alpha\overline{u}\in A$. If $\alpha=1$ then $\overline{u}=\alpha\overline
{u}\in A$. If $\alpha\in(0,1)$ then $\overline{u}=\alpha\overline{u}%
+(1-\alpha)\overline{u}\in A+\operatorname*{int}K\subset A$. Therefore,
$\operatorname*{int}K\cap\operatorname*{conv}A\subset A$, whence
$\operatorname*{conv}A=K\cap\operatorname*{conv}A=\operatorname*{cl}\left(
\operatorname*{int}K\cap\operatorname*{conv}A\right)  \subset
\operatorname*{cl}A=A$. Therefore, $A$ is convex.

Assume that (\ref{r-sa}) does not hold. Then there exist $a,a^{\prime}\in
A\cap\operatorname*{int}K$ with $a\neq a^{\prime}$ and$\ \lambda\in(0,1)$ such
that $\overline{a}:=\lambda a+(1-\lambda)a^{\prime}\in\operatorname*{bd}A$; of
course, $\overline{a}\in\operatorname*{int}K$. By (\ref{fp-ssc}) we have that
$\overline{a}+(K\setminus\{0\})\subset\operatorname*{int}K\cap
\operatorname*{int}_{K}A=\operatorname*{int}A$. Because $A$ is convex, as
above (with $\overline{a}$ instead of $\overline{u}$), there exists $x^{\ast
}\in-K^{\#}=\operatorname*{int}K^{-}$ such that $\nabla\sigma_{A}(x^{\ast
})=\overline{a}$. We have that
\[
\sigma_{A}(x^{\ast})=\left\langle \overline{a},x^{\ast}\right\rangle
=\lambda\left\langle a,x^{\ast}\right\rangle +(1-\lambda)\left\langle
a^{\prime},x^{\ast}\right\rangle \leq\lambda\sigma_{A}(x^{\ast})+(1-\lambda
)\sigma_{A}(x^{\ast})=\sigma_{A}(x^{\ast}),
\]
whence $\left\langle a,x^{\ast}\right\rangle =\left\langle a^{\prime},x^{\ast
}\right\rangle =\sigma_{A}(x^{\ast})$. Hence $a,a^{\prime}\in\partial
\sigma_{A}(x^{\ast})=\{\overline{a}\}$ which yields the contradiction
$a=a^{\prime}$. Therefore, (\ref{r-sa}) holds.

Assume now that $\dim X=2$ and $(\ref{r-sa})$ holds. Hence $K=\mathbb{R}%
_{+}x_{1}+\mathbb{R}_{+}x_{2}$ with $x_{1},x_{2}\in X$ linearly independent.
From $(\ref{r-sa})$ we have obviously that $A\cap\operatorname*{int}K$ is
convex, and so, from the last equality in (\ref{r-intA}) we get the convexity
of $A.$

Let $x,x^{\prime}\in A$ with $x\neq x^{\prime}$ and $\lambda\in(0,1)$. Set
$x^{\prime\prime}:=\lambda x+(1-\lambda)x^{\prime}$. Assume that $x^{\prime
}\in\operatorname*{int}K$; then $u:=\tfrac{1}{2}x+\tfrac{1}{2}x^{\prime\prime
}=\frac{1+\lambda}{2}x+\frac{1-\lambda}{2}x^{\prime}\in A\cap
\operatorname*{int}K$, and so $x^{\prime\prime}=\frac{2\lambda}{1+\lambda
}u+\frac{1-\lambda}{1+\lambda}x^{\prime}\in\operatorname*{int}%
A=A+\operatorname*{int}K\subset A+(K\setminus\{0\})=A+(A_{\infty}%
\setminus\{0\}).$

Assume now that $x,x^{\prime}\in\operatorname*{bd}K=\mathbb{R}_{+}x_{1}%
\cup\mathbb{R}_{+}x_{2}$. If $x,x^{\prime}\in\mathbb{R}_{+}x_{1}$ then
$x=\alpha x_{1}$, $x^{\prime}=\alpha^{\prime}x_{1}$ with $\alpha
,\alpha^{\prime}>0$. Letting $\alpha>\alpha^{\prime}$ we have that
$x^{\prime\prime}=x^{\prime}+\lambda(\alpha-\alpha^{\prime})x_{1}\in
A+(K\setminus\{0\})$. If $x\in\mathbb{R}_{+}x_{1}$ and $x^{\prime}%
\in\mathbb{R}_{+}x_{2}$ then $x=\alpha x_{1}$, $x^{\prime}=\alpha^{\prime
}x_{2}$ with $\alpha,\alpha^{\prime}>0$. Take $u:=\frac{1+\lambda}{2}%
x+\frac{1-\lambda}{2}x^{\prime}\in A\cap\left(  \mathbb{P}x_{1}+\mathbb{P}%
x_{2}\right)  $ and $v:=\frac{\lambda}{2}x+\frac{2-\lambda}{2}x^{\prime}\in
A\cap\left(  \mathbb{P}x_{1}+\mathbb{P}x_{2}\right)  $. Since
$\operatorname*{int}K=\mathbb{P}x_{1}+\mathbb{P}x_{2}$, we obtain that
$x^{\prime\prime}=\lambda u+(1-\lambda)v\in\operatorname*{int}A\subset
A+(A_{\infty}\setminus\{0\})$. Hence $(\ref{rs2})$ is verified. Applying
Proposition \ref{ps2} we obtain that $\sigma_{A}$ is differentiable on
$\operatorname*{int}K^{-}$. The proof is complete. \hfill$\square$

\medskip

Note that $A:=x_{0}+K$ with $x_{0}\in K\setminus\{0\}$ verifies condition (H)
and $\sigma_{A}$ is differentiable on $\operatorname*{int}K^{-}$. However,
$(\ref{fp-ssc})$ is not verified. So, condition $(\ref{fp-ssc})$ is far from
being necessary for the differentiability of $\sigma_{A}.$

(Assume that $x_{0}\in K\setminus\{0\}$ and $A:=x_{0}+K$ verifies
$(\ref{fp-ssc})$. If $x_{0}\in\operatorname*{int}K$, then $A\subset
\operatorname*{int}K$, and so, by (\ref{r-intA}), $\operatorname*{int}%
A=\operatorname*{int}_{K}A$; thus $x_{0}+(K\setminus\{0\})\subset
\operatorname*{int}A=x_{0}+\operatorname*{int}K$, whence $K\setminus
\{0\}\subset\operatorname*{int}K$, a contradiction. Hence $x_{0}%
\in\operatorname*{bd}K$. Then there exists $u\in\operatorname*{int}K$ and
$t>0$ such that $(1+t)u-tx_{0}\notin K$. Otherwise $\operatorname*{int}%
K-x_{0}\subset K$, whence $K\subset x_{0}+K$. This implies the
contradiction $-x_{0}\in K$. It follows that there exist
$u\in\operatorname*{int}K$ and $t>0$ such that
$u_{0}:=(1+t)u-tx_{0}\in\operatorname*{bd}K$. Clearly, $u_{0}\neq0$;
else $x_{0}=(1+t^{-1})u\in\operatorname*{int}K$, a contradiction. It
follows that $x_{0}+t^{-1}u_{0}=(1+t^{-1})u\in
\operatorname*{int}K\cap\operatorname*{int}_{K}A=\operatorname*{int}%
A=x_{0}+\operatorname*{int}K$, and so $u_{0}\in\operatorname*{int}K$, a contradiction.)

A condition which is slightly stronger than (\ref{r-sa}) is sufficient for the
differentiability of $\sigma_{A}.$

\begin{proposition}
\label{prop-suf}Assume that
\begin{equation}
\forall a,a^{\prime}\in A,\ a\neq a^{\prime},\ \forall\lambda\in(0,1):\lambda
a+(1-\lambda)a^{\prime}\in\operatorname*{int}\nolimits_{K}A. \label{r-sas}%
\end{equation}
Then $A$ is convex, $\sigma_{A}$ is differentiable on $\operatorname*{int}%
K^{-}$ and $(\ref{fp-ssc})$ holds.
\end{proposition}

Proof. From (\ref{r-sas}) and $\operatorname*{int}\nolimits_{K}A\subset A$ we
get the convexity of $A$. Let now $a\in A$ and $k\in K\setminus\{0\}$. Then
$a\neq a+2k\in A$. From (\ref{r-sas}) we obtain that $a+k=\tfrac{1}{2}%
a+\tfrac{1}{2}(a+2k)\in\operatorname*{int}\nolimits_{K}A$, and so
(\ref{fp-ssc}) holds.

Assume that $\sigma_{A}$ is not differentiable on $\operatorname*{int}K^{-}$.
Then there exist $\overline{x}^{\ast}\in\operatorname*{int}K^{-}$ and
$a,a^{\prime}\in\partial\sigma_{A}(\overline{x}^{\ast})\subset
\operatorname*{bd}A\subset A$ such that $a\neq a^{\prime}$. By (\ref{r-sas})
we have that $0\neq\overline{a}:=\tfrac{1}{2}a+\tfrac{1}{2}a^{\prime}%
\in\partial\sigma_{A}(\overline{x}^{\ast})\cap\operatorname*{int}%
\nolimits_{K}A$. In particular, $\left\langle \overline{a},\overline{x}^{\ast
}\right\rangle =\sigma_{A}(\overline{x}^{\ast})$. Then there exists $\alpha
\in(0,1)$ such that $\alpha\overline{a}\in A\subset K\setminus\{0\}$, and so
$\left\langle \overline{a},\overline{x}^{\ast}\right\rangle <0$. This is in
contradiction with $\left\langle \overline{a},\overline{x}^{\ast}\right\rangle
=\sigma_{A}(\overline{x}^{\ast})\geq\left\langle \alpha\overline{a}%
,\overline{x}^{\ast}\right\rangle $. Therefore, $\sigma_{A}$ is differentiable
on $\operatorname*{int}K^{-}$. \hfill$\square$

\medskip

The example above (that is $A:=x_{0}+K$ with $x_{0}\in K\setminus\{0\}$) shows
that the condition (\ref{r-sas}) is not necessary for the differentiability of
$\sigma_{A}$ on $\operatorname*{int}K^{-}.$

Applying Theorem \ref{Fact 7}, we obtain that $(\ref{r-sas})\Rightarrow
(\ref{fps})$; the advantage of $(\ref{r-sas})$ is that this condition is more
intuitive (and quite easy to be verified).

As seen above, $\operatorname*{int}A=\operatorname*{int}_{K}A$ and
$\operatorname*{bd}A=\operatorname*{bd}_{K}A$ when $A\subset
\operatorname*{int}K$; in this case condition (\ref{fp-ssc}) becomes
\begin{equation}
A+(K\setminus\{0\})\subset\operatorname*{int}A. \label{r-s}%
\end{equation}
Conversely, if $A$ verifies (\ref{r-s}), then $A\subset\operatorname*{int}K$.
Indeed, in the contrary case there exists $k\in A\setminus\operatorname*{int}%
K\subset K\setminus\{0\}$. From (\ref{r-s}) we get $2k=k+k\in
\operatorname*{int}A\subset\operatorname*{int}K$, whence the contradiction
$k\in\operatorname*{int}K.$

The next result characterizes the differentiability of $\sigma_{A}$ for
$A\subset\operatorname*{int}K.$

\begin{corollary}
\label{cor11}Assume that $A\subset\operatorname*{int}K$ and $(\ref{r-s})$
holds. Then $\sigma_{A}$ is differentiable on $\operatorname*{int}K^{-}$ if
and only if $A$ is convex and
\begin{equation}
\forall a,a^{\prime}\in\operatorname*{bd}A,\ a\neq a^{\prime},\ \forall
\lambda\in(0,1):\lambda a+(1-\lambda)a^{\prime}\notin\operatorname*{bd}A.
\label{r-sb}%
\end{equation}

\end{corollary}

Proof. From the observation above we have that $\operatorname*{int}%
A=\operatorname*{int}_{K}A$, and so (\ref{fp-ssc}) holds. Assume first that
$\sigma_{A}$ is differentiable on $\operatorname*{int}K^{-}$. By Theorem
\ref{Fact 14} we have that $A$ is convex and (\ref{r-sa}) holds. Take
$a,a^{\prime}\in\operatorname*{bd}A\subset\operatorname*{int}K$ with $a\neq
a^{\prime}$ and$\ \lambda\in(0,1)$. By (\ref{r-sa}) we have that $\overline
{a}:=\lambda a+(1-\lambda)a^{\prime}\in\operatorname*{int}A$. Therefore,
$\overline{a}\notin\operatorname*{bd}A$, and so (\ref{r-sb}) holds.

Assume now that $A$ is convex and (\ref{r-sb}) holds. Then $A$ is strictly
convex, and so, using Theorem \ref{ps1}, $\sigma_{A}$ is differentiable on
$\operatorname*{dom}\partial\sigma_{A}\setminus(\operatorname*{lin}%
_{0}A)^{\perp}=\operatorname*{dom}\partial\sigma_{A}\setminus\{0\}\supset
\operatorname*{int}K^{-}$. \hfill$\square$

\medskip

Note that for obtaining Proposition \ref{prop-suf}, or relation (\ref{r-sa})
when $A$ is convex and $\sigma_{A}$ is differentiable on $\operatorname*{int}%
K^{-}$ $(=\operatorname*{int}(\operatorname*{dom}\sigma_{A}))$ in Theorem
\ref{Fact 14}, it is not possible to use anyone of the results in Section 3.

In the sequel we use the convention $0A:=A_{\infty}=K$. Then%
\begin{equation}
\lbrack\alpha,\infty)A=\alpha A\quad\forall\alpha\in\mathbb{R}_{+}.
\label{r11}%
\end{equation}
Indeed, the equality is obvious for $\alpha=0$ (since $A\subset K$ and
$0A=K$). Let $\alpha>0$ and take $t\in\lbrack\alpha,\infty)$, $x\in A$; then
$tx=\alpha\left[  x+(t/\alpha-1)x\right]  \in\alpha\left(  A+K\right)  =\alpha
A.$

The set $\left\{  t\geq0\mid x\in tA\right\}  $ is a compact interval
containing $0$. Using the facts that $0A=K$, $0\notin A=\operatorname*{cl}A,$
and (\ref{r11}), we obtain that the function%
\begin{equation}
F_{A}:K\rightarrow\mathbb{R}_{+},\quad F_{A}(x):=\max\left\{  t\geq0\mid x\in
tA\right\}  . \label{r12}%
\end{equation}
is well defined. The function $F_{A}$ is the restriction to $K$ of the
function $\beta_{A}$ considered in \cite{PenZal:00} in a more general setting.

In the following proposition we mention some properties of $F_{A}.$

\begin{theorem}
\label{prop-fa}Let $A$ and $K$ be as above.

\emph{(i)} $F_{A}(tx)=tF_{A}(x)$ for all $x\in K$ and $t\geq0$, $x\in
F_{A}(x)\cdot A$ for every $x\in K$, and
\begin{equation}
\operatorname*{int}K\subset\mathbb{P}A=\{x\in K\mid F_{A}(x)>0\}. \label{r14}%
\end{equation}

\emph{(ii)} One has%
\begin{equation}
\{x\in K\mid F_{A}(x)\geq\gamma\}=\gamma A\quad\forall\gamma\in\mathbb{R}_{+}.
\label{r13}%
\end{equation}
Consequently, $F_{A}$ is upper semi\-continuous (usc for short).

\emph{(iii)} $A$ is convex iff $F_{A}$ is quasiconcave iff $F_{A}$ is concave.

\emph{(iv)} $F_{A}(x^{\prime})\geq F_{A}(x)$ for all $x,x^{\prime}\in K$ with
$x^{\prime}\geqq_{K}x.$

\emph{(v) }$F_{A}$ is continuous on $(K\setminus\mathbb{P}A)\cup
\operatorname*{int}K$. Therefore, $F_{A}$ is continuous whenever
$A\subset\operatorname*{int}K$.

\emph{(vi)} If $K$ is polyhedral, then $F_{A}$ is continuous.

\emph{(vii)} \emph{a)} $(\ref{r-s})$ holds iff $x^{\prime}\geq_{K}%
x\in\operatorname*{int}K$ implies $F_{A}(x^{\prime})>F_{A}(x)$. \emph{b)}
$(\ref{fp-ssc})$ holds iff $F_{A}$ is continuous and $x^{\prime}\geq_{K}%
x\in\mathbb{P}A$ implies $F_{A}(x^{\prime})>F_{A}(x).$

\emph{(viii)} \emph{a)} $(\ref{r-sa})$ holds iff $F_{A}$ is strictly
quasi-concave on $\operatorname*{int}K$. \emph{b)} If $(\ref{r-sas})$ holds,
then $F_{A}$ is strictly quasi-concave on $\mathbb{P}A$; conversely, if
$F_{A}$ is continuous and strictly quasi-concave on $\mathbb{P}A$ then
$(\ref{r-sas})$ holds.
\end{theorem}

Proof. (i) Since $0\in0A\setminus\mathbb{P}A$ we have that $F_{A}(0)=0$. The
relation $F_{A}(tx)=tF_{A}(x)$ for $x\in K$ and $t>0$ follows immediately from
the definition. Also the relation $x\in F_{A}(x)\cdot A$ follows from the very
definition of $F_{A}$ (the supremum being attained).

The equality in (\ref{r14}) is obvious. Assume that there exists
$x\in(\operatorname*{int}K)\setminus\mathbb{P}A$. Then $\mathbb{P}x\cap
A=\emptyset$, and so $\mathbb{P}x\cap(a+K)=\emptyset$, where $a\in A$ is a
fixed element. Since $\mathbb{P}x$ and $a+K$ are convex sets, there exists
$x^{\ast}\in X^{\ast}\setminus\{0\}$ such that $\left\langle tx,x^{\ast
}\right\rangle \leq\left\langle a+u,x^{\ast}\right\rangle $ for all $t>0$ and
$u\in K$. It follows that $x^{\ast}\in K^{+}$ and $\left\langle x,x^{\ast
}\right\rangle \leq0$. Since $x\in\operatorname*{int}K$ we get the
contradiction $x^{\ast}=0.$

(ii) The inclusion $\supset$ in (\ref{r13}) is obvious. The converse
inclusion is immediate from (\ref{r11}) and the fact that $x\in
F_{A}(x)\cdot A$ for every $x\in K$. Because $\gamma A$ is closed
for every $\gamma\geq0$, from (\ref{r13}) we obtain that $F_{A}$ is
usc.

(iii) The first equivalence follows from (\ref{r13}). Moreover, if
$F_{A}$ is concave, clearly $F_{A}$ is quasiconcave. Assume that
$F_{A}$ is quasiconcave. Consider
$f:X\rightarrow\overline{\mathbb{R}}$ defined by $f(x):=-F_{A}(x)$
for $x\in K$ and $f(x):=+\infty$ for $x\in X\setminus K$. Then $f$
is quasiconvex. Because $F_{A}$ is usc and $\operatorname*{dom}f=K$
is closed we have that $f$ is lsc. Moreover, from (i) we have that
$f(tx)=tf(x)$ for all $t\in\mathbb{P}$ and $x\in X$, and $\{x\in
X\mid f(x)<0\}=\{x\in K\mid
F_{A}(x)>0\}\supset\operatorname*{int}K$, whence $\operatorname*{dom}%
f=K=\operatorname*{cl}\{x\in X\mid f(x)<0\}$. Applying \cite[Thm.\ 2.2.2]%
{Zal:02} we obtain that $f$ is sublinear; in particular, $f$ is convex, and so
$F_{A}$ is concave.

(iv) Take $x,x^{\prime}\in K$ with $x^{\prime}\geqq_{K}x$. If
$\gamma :=F_{A}(x)=0$ then clearly $F_{A}(x^{\prime})\geq\gamma$.
Else, $\gamma>0$ and $x\in\gamma A$; hence $x^{\prime}\in
x+K\subset\gamma A+K=\gamma(A+K)=\gamma A$, and so
$F_{A}(x^{\prime})\geq\gamma$ by (ii).

(v) Let $x\in\operatorname*{int}K$; then $\gamma:=F_{A}(x)>0$. Take
$0<\mu<\gamma$. Then
\[
\mu^{-1}x=\gamma^{-1}x+(\mu^{-1}-\gamma^{-1})x\in A+\operatorname*{int}%
K\subset\operatorname*{int}A.
\]
It follows that $A$ is a neighborhood of $\mu^{-1}x$, whence $V:=\mu A$ is a
neighborhood of $x$. Since $F_{A}(x^{\prime})\geq\mu$ for every $x^{\prime}\in
V$, we have that $F_{A}$ is lsc at $x$. By (ii) we get the continuity of
$F_{A}$ at $x.$

Take $x\in K\setminus\mathbb{P}A$; from (\ref{r14}) we have that
$F_{A}(x)=0=\inf F_{A}$, and so $F_{A}$ is lsc at $x$. Since $F_{A}$ is usc,
we have that $F_{A}$ is continuous at $x$. Hence $F_{A}$ is continuous at any
$x\in(K\setminus\mathbb{P}A)\cup\operatorname*{int}K.$

If $A\subset\operatorname*{int}K$, then $\mathbb{P}A=\operatorname*{int}K$,
and so $(K\setminus\mathbb{P}A)\cup\operatorname*{int}K=K.$

(vi) There exists $(x_{i}^{\ast})_{i\in\overline{1,m}}\subset X^{\ast
}\setminus\{0\}$ such that $K=\{x\in X\mid\left\langle x,x_{i}^{\ast
}\right\rangle \geq0\ \forall i\in\overline{1,m}\}$, $\operatorname*{int}%
K=\{x\in X\mid\left\langle x,x_{i}^{\ast}\right\rangle >0\ \forall
i\in\overline{1,m}\}\neq\emptyset$ and $\cap_{i=1}^{m}\ker x_{i}^{\ast
}=\{0\}.$

Take $a\in\mathbb{P}A$ and set $I:=\{i\in\overline{1,m}\mid\left\langle
a,x_{i}^{\ast}\right\rangle >0\}$ $(\neq\emptyset)$. Let $\gamma:=F_{A}(a)$
$(>0)$ and take $\mu\in(0,\gamma)$; clearly $a\in\gamma A$. There exists a
neighborhood $V$ of $a$ such that $\left\langle x,x_{i}^{\ast}\right\rangle
\geq\mu\gamma^{-1}\left\langle a,x_{i}^{\ast}\right\rangle $ for all $x\in V$
and $i\in I$. Then for each $x\in K\cap V$ and each $i\in\overline{1,m}$ we
have that $\left\langle \mu^{-1}\gamma x-a,x_{i}^{\ast}\right\rangle \geq0.$
Thus $\mu^{-1}x\in\gamma^{-1}a+K\subset A$. Hence $F_{A}(x)\geq\mu$ for every
$x\in K\cap V$, and so $F_{A}$ is lsc at $x$. It follows that $F_{A}$ is
continuous at $x.$

(vii) b) \textquotedblleft$\Longrightarrow$\textquotedblright\ If $x\in
K\setminus\mathbb{P}A$ then $F_{A}(x)=0=\inf F_{A}$, and so $F_{A}$ is lsc
(hence continuous by (ii)) at $x$. Let $x\in\mathbb{P}A$ $(\subset
K\setminus\{0\})$ and take $\gamma:=F_{A}(x)>0$. Consider $0<\mu<\gamma$.
Then, as in (v), we get $\mu^{-1}x\in A+K\setminus\{0\}\subset
\operatorname*{int}\nolimits_{K}A$. Hence $V:=\mu A$ is a neighborhood (in
$K$) of $x$. Since $F_{A}(x^{\prime})\geq\mu$ for every $x^{\prime}\in V$, we
obtain that $F_{A}$ is lsc at $x.$

Let now $x^{\prime}\geq_{K}x\in\mathbb{P}A$ and take $\gamma:=F_{A}(x)>0$;
then $x^{\prime}=x+k$ $(\in K\setminus\{0\})$ for some $k\in K\setminus\{0\}.$
It follows that $\gamma^{-1}x^{\prime}=\gamma^{-1}x+\gamma^{-1}k\in
\operatorname*{int}\nolimits_{K}A$, and so there exists a neighborhood $V$ of
$\gamma^{-1}x^{\prime}$ such that $K\cap V\subset A$. Then there exists
$\mu>\gamma$ such that $\mu^{-1}x^{\prime}\in K\cap V\subset A$, whence
$F_{A}(x^{\prime})\geq\mu>\gamma.$

\textquotedblleft$\Longleftarrow$\textquotedblright\ Consider $x\in A$, $k\in
K\setminus\{0\}$ and $x^{\prime}:=x+k$. Then $x^{\prime}\geq_{K}x\in
\mathbb{P}A$, and so $F_{A}(x^{\prime})>F_{A}(x)\geq1$. Because $F_{A}$ is
continuous at $x^{\prime}$, there exists a neighborhood $V$ of $x^{\prime}$
such that $F_{A}(x^{\prime\prime})\geq1$ for every $x^{\prime\prime}\in K\cap
V$. It follows that $K\cap V\subset A$, and so $x^{\prime}\in
\operatorname*{int}\nolimits_{K}A.$

The proof of a) is similar.

(viii) b) Assume that $(\ref{r-sas})$ holds. Take $x,x^{\prime}\in\mathbb{P}A$
with $x\neq x^{\prime}$ and $\lambda\in(0,1)$. We may (and do) assume that
$F_{A}(x^{\prime})\geq F_{A}(x)=:\gamma>0$. Then $x^{\prime},x\in\gamma A$. It
follows that $\gamma^{-1}x^{\prime\prime}\in\operatorname*{int}\nolimits_{K}%
A$, where $x^{\prime\prime}:=(\lambda x+(1-\lambda)x^{\prime})$. As in the
proof of (vii) b) above, there exists $\mu>\gamma$ such that $\mu
^{-1}x^{\prime\prime}\in A$, whence $F_{A}(x^{\prime\prime})\geq\mu
>\gamma=F_{A}(x)$. Hence $F_{A}$ is strictly quasi-concave on $\mathbb{P}A.$

Assume now that $F_{A}$ is continuous and strictly quasi-concave on
$\mathbb{P}A$. Take $x,x^{\prime}\in A$ with $x\neq x^{\prime}$ and
$\lambda\in(0,1)$. Then $F_{A}(x)$, $F_{A}(x^{\prime})\geq1$.
Because $F_{A}$ is strictly quasi-concave, we have that
$F_{A}(x^{\prime\prime})>1$, where $x^{\prime\prime}:=\lambda
x+(1-\lambda)x^{\prime}$. From (\ref{r14}) we have that
$x^{\prime\prime}\in\mathbb{P}A$, and so $F_{A}$ is continuous at
$x^{\prime\prime}$. Then there exists a neighborhood $V$ of
$x^{\prime\prime}$ such that $F_{A}(y)\geq1$ for every $y\in K\cap
V$, whence $K\cap V\subset A$. This shows that
$x^{\prime\prime}\in\operatorname*{int}\nolimits_{K}A.$

The proof of a) is similar; take into account that $F_{A}$ is continuous on
$\operatorname*{int}K$. \hfill$\square$

\medskip

The next examples show that in several results the converse implications are
not valid.

\begin{example}
\label{ex3}\emph{(a) The set $A:=a+K$ with $a\in K\setminus\{0\}$
does not verify condition $(\ref{fp-ssc})$, but $\sigma_{A}$ is
differentiable on $\operatorname*{int}K^{-}$. Also the condition
$(\ref{r-s})$ is not necessary for the differentiability of
$\sigma_{A}$ when $A\subset\operatorname*{int}K$ (take
$a\in\operatorname*{int}K$).}

\emph{(b) The set $A:=\{(x_{1},x_{2})\in\mathbb{R}_{+}^{2}\mid
x_{1}+x_{2}\geq1\}$ verifies condition $(\ref{fp-ssc})$ for
$K:=\mathbb{R}_{+}^{2}$, but
$\sigma_{A}$ is not differentiable on $\operatorname*{int}K^{-}=-\mathbb{R}%
_{++}^{2}.$}

\emph{(c) Let $A:=a+K$ with $K:=\{(x_{1},x_{2},x_{3})\in\mathbb{R}^{3}\mid x_{3}%
\geq\sqrt{(x_{1})^{2}+(x_{2})^{2}}\}$ and $a\in(\operatorname*{bd}%
K)\setminus\{0\}$. Then
$\mathbb{P}A=\mathbb{P}a\cup\operatorname*{int}K$ and $F_{A}$ is not
continuous at $x\in K$ iff $x\in\mathbb{P}a$. (This fact shows that
the polyhedrality of $K$ in Theorem \ref{prop-fa}~(vi) is
essential.) Moreover, $F_{A}$ is sublinear.}

\emph{(d) Let $A:=a+K$ with $K:=\{x\in X\mid\left\langle
x,x_{i}^{\ast}\right\rangle \geq0\ \forall i\in\overline{1,m}\}$,
where $(x_{i}^{\ast})_{i\in
\overline{1,m}}\subset X^{\ast}\setminus\{0\}$ are such that $\cap_{i=1}%
^{m}\ker x_{i}^{\ast}=\{0\}$ and $(\operatorname*{int}K=)$ $\{x\in
X\mid\left\langle x,x_{i}^{\ast}\right\rangle >0\ \forall
i\in\overline {1,m}\}\neq\emptyset$, and $a\in K\setminus\{0\}$.
Then $F_{A}(x)=\min\left\{ \left\langle x,x_{i}^{\ast}\right\rangle
/\left\langle a,x_{i}^{\ast }\right\rangle \mid\left\langle
a,x_{i}^{\ast}\right\rangle >0\right\}  $ for $x\in K$. For
$K:=\mathbb{R}_{++}^{p}$, $F_{A}$ is the Leontieff production
function (and $F_{A}(x)=\min\left\{  x_{i}/a_{i}\mid a_{i}>0\right\}
$).}
\end{example}

(a) This example was considered before Proposition \ref{prop-suf}.
Moreover,
$\sigma_{A}(x^{\ast})=\left\langle a,x^{\ast}\right\rangle +\iota_{K^{-}%
}(x^{\ast})$.

(b) In this case $A+\left(  K\setminus\{0\}\right)  =\{(x_{1},x_{2}%
)\in\mathbb{R}_{+}^{2}\mid x_{1}+x_{2}>1\}=\operatorname*{int}\nolimits_{K}A$
and $\sigma_{A}(u_{1},u_{2})=\max\{u_{1},u_{2}\}+\iota_{-\mathbb{R}_{+}^{2}%
}(u_{1},u_{2}).$

(c) Let
$a:=(a_{1},a_{2},a_{3})\in(\operatorname*{bd}K)\setminus\{0\}$; then
$a_{3}=\sqrt{(a_{1})^{2}+(a_{2})^{2}}>0$. Then for every $x:=(x_{1}%
,x_{2},x_{3})\in K$ we have that $a_{1}x_{1}+a_{2}x_{2}\leq\sqrt{(a_{1}%
)^{2}+(a_{2})^{2}}\sqrt{(x_{1})^{2}+(x_{2})^{2}}\leq a_{3}x_{3}$,
with equality iff $x\in\mathbb{P}a$. After some computation we get
\[
F_{A}(x)=\left\{
\begin{array}
[c]{ll}%
0 & \text{if }x\in(\operatorname*{bd}K)\setminus\mathbb{P}a,\\
x_{3}/a_{3} & \text{if }x\in\mathbb{P}a,\\
\frac{(x_{3})^{2}-(x_{1})^{2}-(x_{2})^{2}}{2\left(  a_{3}x_{3}-a_{1}%
x_{1}-a_{2}x_{2}\right)  } & \text{if }x\in\operatorname*{int}K.
\end{array}
\right.
\]
Clearly, $F_{A}$ is continuous at $x\in K$ iff $x\in K\setminus\mathbb{P}a.$

(d) Because $a\in K\setminus\{0\}$, the set $I:=\{i\in\overline{1,m}%
\mid\left\langle a,x_{i}^{\ast}\right\rangle >0\}\neq\emptyset$. Let
$x\in K$ and $t>0$; then $x\in t(a+K)$ iff $x\geq ta$ iff
$t\leq\min\left\{ \left\langle x,x_{i}^{\ast}\right\rangle
/\left\langle a,x_{i}^{\ast }\right\rangle \mid i\in I\right\}  $.
Hence $F_{A}(x)=\min\left\{ \left\langle x,x_{i}^{\ast}\right\rangle
/\left\langle a,x_{i}^{\ast }\right\rangle \mid i\in I\right\}  $
for every $x\in K.$

\bigskip

Taking into account Theorem \ref{prop-fa}~(vi) and Example \ref{ex3}
(c) one can ask if for any non polyhedral (pointed closed convex)
cone $K$ one can find $A\subset K$ such that $(A,K)$ verify (H) and
$F_{A}$ be not continuous.

In the sequel we give a positive answer. We begin with the next
auxiliary result; the results from Convex analysis used in its proof
can be found in \cite{Roc:70} or \cite{Zal:02}.

\begin{lemma}
\label{lem2}Let $C\subset X$ be a closed convex set and $x_{0}\in C$
be such that $\mathbb{R}_{+}(C-x_{0})$ is not closed. Then there
exists a sequence $(x_{n})_{n\geq1}$ converging to $x_{0}$ such that
$x_{n}\in C\setminus\left[  (1-\lambda)x_{0}+\lambda C\right]  $ for
all $\lambda \in(0,1)$ and $n\ge 1$.
\end{lemma}

Proof. We may (and do) assume that $x_{0}=0$, $\dim\left(  \operatorname*{lin}%
(C-x_{0})\right)  \geq2$ and $\operatorname*{int}C\neq\emptyset$.
Indeed, if
$x_{0}\neq0$ we replace $C$ by $C-x_{0}$; if $\dim\left(  \operatorname*{lin}%
(C-x_{0})\right)  <2$ then $\mathbb{R}_{+}(C-x_{0})=\mathbb{R}_{+}C$ is
closed; if $\operatorname*{int}C=\emptyset$ we replace $X$ by
$\operatorname*{lin}C$.

Consider $u\in\left[  \operatorname*{cl}(\mathbb{R}_{+}C)\right]
\setminus(\mathbb{R}_{+}C)$ with $\left\Vert u\right\Vert =1$ and
$a\in\operatorname*{int}C$; it follows that
$0\in\operatorname*{bd}C$ and $sa\in\operatorname*{int}C$ for every
$s\in(0,1]$. Clearly, $u$ and $a$ are linearly independent. There
exists a basis $\{e_{1},e_{2},\ldots,e_{p}\}$ of
$X$ such that $e_{1}=u$ and $e_{p}=a$. Set $Y:=\operatorname*{lin}%
\{e_{1},e_{2},\ldots,e_{p-1}\}$ and endow $Y$ with the induced norm; then
$X=Y\oplus\mathbb{R}e_{p}=Y\oplus\mathbb{R}a$ and $T:Y\times\mathbb{R}%
\rightarrow X$, $T(y,s):=y+sa$ is an isomorphism of (normed) linear
spaces. Consider the function
$\varphi:Y\rightarrow\overline{\mathbb{R}}$ defined by
$\varphi(y):=\inf\{s\in\mathbb{R}\mid y+sa\in C\}$. Then $\varphi$
is convex, $\varphi(0)=0$ and $C\subset
T(\operatorname*{epi}\varphi)$. Since
$T(0,1)=a\in\operatorname*{int}C$, it follows that $(0,1)\in
\operatorname*{int}(\operatorname*{epi}\varphi)$, and so $\varphi$
is proper and continuous at $0$; therefore, $\varphi$ is Lipschitz
on a neighborhood of $0$. Because $C$ is closed, it follows that
$y+\varphi(y)a\in C$ for every $y\in\operatorname*{dom}\varphi.$
Since $T(0,1)=a\in\operatorname*{int}C$, there exists
$\varepsilon_{0}>0$ such that $T(y,1)\in C$ for every $y\in Y$ with
$\left\Vert y\right\Vert \leq\varepsilon_{0}.$

(a) For $t\in(0,\varepsilon_{0}]$ we have that $s:=\varphi(tu)>0$.
Indeed, in the contrary case $s\leq0$, and so, for
$\lambda:=s/(s-1)\in\lbrack0,1)$, we get the contradiction
$tu=(1-\lambda)(tu+sa)+\lambda(tu+a)\in C$.

(b) We have that  $\varphi^{\prime}(0,u)=0$. Indeed, because $u\in
\operatorname*{cl}(\mathbb{R}_{+}C)$ (and $C$ is convex with $0\in
C$), there exist the sequences $(t_{n})\subset\mathbb{P}$ and
$(a_{n})\subset C$ such that $a_{n}:=y_{n}+s_{n}a\rightarrow0$ and
$t_{n}^{-1}a_{n}\rightarrow u$. Hence $y_{n}\rightarrow0,$
$t_{n}\rightarrow0$, $t_{n}^{-1}s_{n}\rightarrow0$ and
$t_{n}^{-1}y_{n}\rightarrow u$. It follows that for $n$ sufficiently
large, $\varphi(t_{n}u)\leq\varphi(y_{n})+L\left\Vert
y_{n}-t_{n}u\right\Vert
\leq s_{n}+L\left\Vert y_{n}-t_{n}u\right\Vert $, and so $0<t_{n}^{-1}%
\varphi(t_{n}u)\leq t_{n}^{-1}s_{n}+L\left\Vert t_{n}^{-1}y_{n}-u\right\Vert $
for $n$ large. Taking the limit we get $0=\lim t_{n}^{-1}\varphi
(t_{n}u)=\varphi^{\prime}(0,u).$

(c) For all $y\in\mathbb{P}u\cap\operatorname*{dom}\varphi$ and
$\lambda \in(0,1)$ we have that $\varphi(\lambda
y)<\lambda\varphi(y)$. In the contrary case there exist such $y$ and
$\lambda$ with $\lambda\varphi(y)\leq \varphi(\lambda
y)=\varphi(\lambda y+(1-\lambda)0)\leq\lambda\varphi
(y)+(1-\lambda)\varphi(0)=\lambda\varphi(y)<\infty$. By the
convexity of $\varphi$ we get $\varphi(\eta y)=\eta\varphi(y)$ for
every $\eta\in(0,1),$
whence $\varphi^{\prime}(0,y)=\varphi(y)$. Since $y=t_{0}u$ with $t_{0}%
\in\mathbb{P}$, from (b), we obtain that
$\varphi(tu)=\varphi(y)=\varphi ^{\prime}(0,t_{0}u)=0$ for every
$t\in\lbrack0,t_{0}]$, contradicting the fact that $\varphi(tu)>0$
for $t\in(0,\varepsilon_{0}]$ (with $\varepsilon_{0}$ from (a)).

Fix now a sequence $(t_{n})\subset(0,\varepsilon_{0}]$ with $t_{n}%
\rightarrow0$, and take $x_{n}:=t_{n}u+\varphi(t_{n}u)a\in C$.
Clearly, $x_{n}\rightarrow0$. Assuming that $x_{n}\in\lambda C$ for
some $\lambda
\in(0,1)$, we obtain that $\varphi(\lambda^{-1}t_{n}u)\leq\lambda^{-1}%
\varphi(t_{n}u)$, and so $\lambda^{-1}t_{n}u\in\mathbb{P}u\cap
\operatorname*{dom}\varphi$. From (c) we get the contradiction
$\varphi
(t_{n}u)<\lambda\varphi(\lambda^{-1}t_{n}u)\leq\varphi(t_{n}u).$
Therefore, $x_{n}\notin\lambda C$ for all $n\geq1$ and
$\lambda\in(0,1)$. \hfill $\square$

\medskip

In the next result we complete Theorem \ref{prop-fa} vi).

\begin{corollary}
\label{cor-cfa}Let $K\subset X$ be a pointed closed convex cone with nonempty
interior. Then $F_{A}$ is continuous for every set $A\subset K$ satisfying
condition (H) if and only if $K$ is polyhedral.
\end{corollary}

Proof. The sufficiency follows from Theorem \ref{prop-fa} vi).

Assume that $K$ is not polyhedral. Using \cite[Prop.\ 2]{ShaNem:03},
there exists $x_{0}\in K$ such that $\mathbb{R}_{+}(K-x_{0})$ is not
closed; of course, $x_{0}\neq0$. Let us take $A:=x_{0}+K$. Clearly,
$F_{A}(x_{0})=1$. By Lemma \ref{lem2}, there exists a sequence
$(x_{n})_{n\geq1}\subset K\setminus\left(
\tfrac{1}{2}x_{0}+\tfrac{1}{2}K\right)  $ with $x_{n}\rightarrow
x_{0}$. Because $x_{n}\notin\tfrac{1}{2}A$, it follows that
$F_{A}(x_{n})\leq\tfrac{1}{2}$ for every $n\geq1$, and so $F_{A}$ is
not continuous at $x_{0}$. \hfill $\square$

\section{Applications to the differentiability of the cost function}

In the literature the problem discussed in the previous section is related to
the cost function associated to a production function $F:\mathbb{R}_{+}%
^{p}\rightarrow\mathbb{R}_{+}$ $(p\geq2)$ satisfying certain
conditions; we denote simply $\geqq$, $\geq$, $>$ the symbols
$\geqq_{\mathbb{R}_{+}^{p}}$, $\geq_{\mathbb{R}_{+}^{p}}$,
$>_{\mathbb{R}_{+}^{p}}$, respectively. Among the properties the
production function $F$ could have we mention first those used in
\cite{FarPri:86}:

\medskip

F.1~$F(0)=0,$

\medskip

F.2~$F(x)\geq F(x^{\prime})$ if $x\geqq x^{\prime},$

\medskip

F.3~$F$ is quasiconcave,

\medskip

F.4~$F$ is upper semi\-continuous.

\medskip

In the economics literature some of the conditions above are strengthened:

\medskip

F.2b~$F(x)>F(x^{\prime})$ if $x\geq x^{\prime},$

\medskip

F.3b~$F$ is strictly quasiconcave,

\medskip

F.4b~$F$ is continuous.

\medskip

The properties of $F_{A}$ (defined in (\ref{r12})) mentioned in Theorem
\ref{prop-fa} suggest the consideration of the following new conditions:

\medskip

F.2c~$F(x)>F(x^{\prime})$ if $x\geq x^{\prime}$ and $F(x^{\prime})>0,$

\medskip

F.2d~$F(x)>F(x^{\prime})$ if $x\geq x^{\prime}\in\mathbb{R}_{++}^{p},$

\medskip

F.3c~$F$ is strictly quasiconcave on $\{x\in\mathbb{R}_{+}^{p}\mid F(x)>0\}$,

\medskip

F.3d~$F$ is strictly quasiconcave on $\mathbb{R}_{++}^{p},$

\medskip

F.4c~$F$ is continuous on $\{x\in\mathbb{R}_{+}^{p}\mid F(x)>0\}$,

\medskip

F.4d~$F$ is continuous on $\mathbb{R}_{++}^{p}$,

\medskip

F.5~$F(x)>0$ if $x\in\mathbb{R}_{++}^{p}$.

\medskip

Note that in the context of the differentiability of the cost functions Sakai
\cite{Sak:73} used conditions F.1, F.2, F.4b and the fact that $F$ is strictly
concave instead of F.3b; Avriel \textit{et.~al.} \cite{AvrDieSchZan:88} used
conditions F.2, F.3, F.4; Saijo \cite{Sai:83} used conditions F.2b and F.4;
Fuchs-Selinger \cite{Fuc:97} used the condition F.2b.

As in \cite{FarPri:86}, set
\[
L(\gamma):=\{x\in\mathbb{R}_{+}^{p}\mid F(x)\geq\gamma\}\quad(\gamma
\in\mathbb{R}_{+}).
\]
It is well known that:

\medskip

F.2 $\iff$ $L(\gamma)+\mathbb{R}_{+}^{p}\subset L(\gamma)$ for every
$\gamma\in\mathbb{R}_{+}$;

\medskip

F.3 $\iff$ $L(\gamma)$ is convex for every $\gamma\in\mathbb{R}_{+}$;

\medskip

F.4 $\iff$ $L(\gamma)$ is closed for every $\gamma\in\mathbb{R}_{+}$.

\medskip

In the next proposition we establish several relations among the conditions
mentioned above.

\begin{proposition}
\label{prop3}Let $F:\mathbb{R}_{+}^{p}\rightarrow\mathbb{R}_{+}.$

\emph{(i)} F.2b $\Rightarrow$ F.2c $\Rightarrow$ F.2, F.2b $\Rightarrow$ F.2d,
(F.2c $\wedge$ F.5) $\Rightarrow$ F.2d, (F.2d $\wedge$ F.4b) $\Rightarrow$ F.2.

\emph{(ii)} F.3b $\Rightarrow$ F.3c $\Rightarrow$ F.3, F.3b $\Rightarrow$
F.3d, (F.3c $\wedge$ F.5) $\Rightarrow$ F.3d, (F.3d $\wedge$ F.4b)
$\Rightarrow$ F.3.

\emph{(iii)} (F.2 $\wedge$ F.3z) $\Rightarrow$ F.2z, z being b, c or d.

\emph{(iv)} F.4b $\Rightarrow$ (F.4c $\wedge$ F.4d $\wedge$ F.4), (F.4c
$\wedge$ F.4) $\Rightarrow$ F.4b, (F.4c $\wedge$ F.5) $\Rightarrow$ F.4d.
\end{proposition}

Proof. (i) Excepting the last one, the implications are obvious.

(F.2d $\wedge$ F.4b) $\Rightarrow$ F.2: Let $x\geq x^{\prime}$; set
$k:=x-x^{\prime}\geq0$. There exists $(x_{n}^{\prime})\subset\mathbb{R}%
_{++}^{p}$ such that $x_{n}^{\prime}\rightarrow x^{\prime}$; then
$x_{n}:=x_{n}^{\prime}+k\in\mathbb{R}_{++}^{p}$ and $x_{n}\rightarrow x.$
Because $x_{n}\geq x_{n}^{\prime}$ we have, by F.2d, that $F(x_{n}%
)>F(x_{n}^{\prime})$ for every $n$, and so $F(x)\geq F(x^{\prime})$ by F.4b.

(ii) Excepting the last one, the implications are obvious.

(F.3d $\wedge$ F.4b) $\Rightarrow$ F.3: Let $x,x^{\prime}\in\mathbb{R}_{+}%
^{p}$ and $\lambda\in(0,1)$. There exist the sequences $(x_{n}),$
$(x_{n}^{\prime})\subset\mathbb{R}_{++}^{p}$ such that
$x_{n}\rightarrow x$ and $x_{n}^{\prime}\rightarrow x^{\prime}$.
Then $F\left(  \lambda x_{n}+(1-\lambda)x_{n}^{\prime}\right)
\geq\min\{F(x_{n}),F(x_{n}^{\prime })\}$ (by F.3d), and so, taking
the limit, we get $F\left(  \lambda x+(1-\lambda)x^{\prime}\right)
\geq\min\{F(x),F(x^{\prime})\}$ (by F.4b).

(iii) (F.2 $\wedge$ F.3z) $\Rightarrow$ F.2z: Let $x^{\prime}\geq x$. Hence
$F(x^{\prime})\geq F(x)$ by $F2$; moreover, $F(x^{\prime})>0$ if $F(x)>0,$
respectively $x^{\prime}\in\mathbb{R}_{++}^{p}$ if $x\in\mathbb{R}_{++}^{p}.$
Then $x\neq x^{\prime\prime}:=2x^{\prime}-x\in\mathbb{R}_{+}^{p}$ and
$x^{\prime}=\tfrac{1}{2}x^{\prime\prime}+\tfrac{1}{2}x$; because
$x^{\prime\prime}\geq x$ we have $F(x^{\prime\prime})\geq F(x)$ by F.2, and
$F(x^{\prime})>F(x)$ by F.3z.

(iv) The implications are obvious. \hfill$\square$

\medskip

When referring to results in the previous sections, in the sequel
$X$ is $\mathbb{R}^{p}$ endowed with the Euclidean norm and
identified with its dual. Because the conditions F.1, F.2 and F.4
seems to be very natural, in the sequel we also assume that $F$
verifies these conditions. In this situation, taking
$K:=\mathbb{R}_{+}^{p}$ and $A:=L(\gamma)$, we have that
$A=\operatorname*{cl}A=A+K\subset K$, and $0\in A$ iff $\gamma=0$.
Hence, if $A\neq\emptyset$ and $\gamma>0$ then $K$ and $A$ verify
condition (H). Set $\Gamma_{F}:=\{\gamma\in\mathbb{P}\mid
L(\gamma)\neq\emptyset\}$; clearly
$\operatorname{Im}F\setminus\{0\}\subset\Gamma_{F}\subset(0,\sup
\operatorname{Im}F]$. If $F$ is continuous then $(0,\sup\operatorname{Im}%
F)\subset\operatorname{Im}F\setminus\{0\}.$

\begin{proposition}
\label{prop2}Let $F:\mathbb{R}_{+}^{p}\rightarrow\mathbb{R}_{+}$ be continuous
and set $K:=\mathbb{R}_{+}^{p}$. Then

\emph{(a)}~F.2c $\Rightarrow$ [$(\ref{fp-ssc})$ with $A:=L(\gamma)$] for every
$\gamma\in\Gamma_{F},$

\emph{(b)}~F.3c $\Rightarrow$ [$(\ref{r-sas})$ with $A:=L(\gamma)$] for every
$\gamma\in\Gamma_{F}.$
\end{proposition}

Proof. (a) Take $\gamma\in\Gamma_{F}$, $x\in L(\gamma)$ and $k\in
K\setminus\{0\}$; hence $F(x)>0$. By F.2b we have that $F(x^{\prime}%
)>F(x)\geq\gamma$, where $x^{\prime}:=x+k$. Because $F$ is continuous at
$x^{\prime}$, there exists a neighborhood $V$ of $x^{\prime}$ such that
$F(u)\geq\gamma$ for every $u\in K\cap V$. Therefore, $K\cap V\subset
L(\gamma)$, which proves that $x^{\prime}\in\operatorname*{int}_{K}L(\gamma)$.
Hence (\ref{fp-ssc}) holds.

(b) Take $\gamma\in\Gamma_{F}$, $x,x^{\prime}\in L(\gamma)$ with $x\neq
x^{\prime}$ and $\lambda\in(0,1)$; assume that $F(x)\geq F(x^{\prime})$. Since
$F$ is strictly quasiconcave on $B:=\{u\mid F(u)>0\}$, we have that $F\left(
x^{\prime\prime}\right)  >F(x^{\prime})\geq\gamma$, where $x^{\prime\prime
}:=\lambda x+(1-\lambda)x^{\prime}$. Since $F$ is continuous on $B$, there
exists a neighborhood $V$ of $x^{\prime}$ such that $F(u)\geq\gamma$ for every
$u\in K\cap V$. Therefore, $K\cap V\subset L(\gamma)$, which proves that
$x^{\prime\prime}\in\operatorname*{int}_{K}L(\gamma)$. Hence (\ref{r-sas})
holds. \hfill$\square$

\medskip

The following questions are quite natural: Are the converse implications in
Proposition \ref{prop2} true? More precisely, if $F:\mathbb{R}_{+}%
^{p}\rightarrow\mathbb{R}_{+}$ is continuous, is it true that F.2c holds if
$L(\gamma)$ satisfies $(\ref{fp-ssc})$ for every $\gamma\in\Gamma_{F}?$ Is it
true that F.3c holds if $L(\gamma)$ satisfies (\ref{r-sas}) for every
$\gamma\in\Gamma_{F}?$

The answer is negative for both questions. For this take
$G:\mathbb{R}_{+}^{p}\rightarrow \mathbb{R}_{+}$ satisfying
conditions F.1, F.2b, F.3b, F.4b and $\sup G>1$; $G$ could be
defined by $G(x_{1},x_{2}):=x_{1}+x_{2}+\sqrt{x_{1}x_{2}}$ for
$(x_{1},x_{2})\in\mathbb{R}_{+}^{2}$. Take also $\varphi:\mathbb{R}%
_{+}\rightarrow\mathbb{R}_{+}$,
$\varphi(t):=\min\{t,\max\{1,t-1\}\}$ and $F:=\varphi\circ G;$
because $\varphi$ is a continuous non decreasing function with
$\varphi(0)=0$, $F$ verifies F.1, F.2b, F.3 and F.4. Then for every
$\gamma\in\mathbb{R}_{+}$ we have that $L_{F}(\gamma)=L_{G}(\gamma)$
for $\gamma\in(0,1]$ and $L_{F}(\gamma)=L_{G}(\gamma+1)$ for
$\gamma\in(1,\infty )$. Hence $L_{F}(\gamma)$ satisfies conditions
$(\ref{fp-ssc})$ and $(\ref{r-sas})$ for every $\gamma\geq0$ with
$L_{F}(\gamma)\neq\emptyset.$
Since $F$ is constant on the nonempty open set $\{(x,y)\in\mathbb{R}_{++}%
^{p}\mid1<G(x,y)<2\}$, we obtain that $F$ satisfies neither F.2c nor F.3c.

The next example shows that we can not replace the continuity of $F$ by its
upper semi\-continuity in Proposition \ref{prop2} (b).

\begin{example}
\label{ex-adsz}\emph{(\cite[Ex.\ 4.4]{AvrDieSchZan:88}) Let $F:\mathbb{R}_{+}%
^{2}\rightarrow\mathbb{R}_{+}$ be defined by%
\[
F(x_{1},x_{2}):=\left\{
\begin{array}
[c]{ll}%
x_{1}x_{2} & \text{if }x_{1}x_{2}<1\text{ or }[\tfrac{1}{2}<x_{1}<2\text{ and
}x_{2}=1/x_{1}],\\
1+x_{1}x_{2} & \text{if }x_{1}x_{2}\geq1\text{ and }x_{1}+x_{2}\geq\tfrac
{5}{2},\\
1+\frac{x_{1}x_{2}-1}{(5/2)x_{1}-1-(x_{1})^{2}} & \text{if }x_{1}x_{2}>1\text{
and }x_{1}+x_{2}<\tfrac{5}{2}.
\end{array}
\right.
\]
The function $F$ is usc, it is strictly quasiconcave on $\mathbb{R}_{++}%
^{2}=\{(x_{1},x_{2})\mid F(x_{1},x_{2})>0\}$, but $F$ is not continuous [for
example, $F$ is not continuous at $(2,\tfrac{1}{2})$]. Hence F.1, F.3c and F.4
hold, but F.4c does not hold. Moreover, F.2c holds by Proposition \ref{prop3}
(iii). However, $L(5/2)=\{(x_{1},x_{2})\in\mathbb{R}_{+}^{2}\mid x_{1}%
x_{2}\geq1,\ x_{1}+x_{2}\geq\tfrac{5}{2}\}$ does not verify
$(\ref{r-sas})$ [or, equivalently, $(\ref{r-sb})$ because
$L(5/2)\subset\mathbb{R}_{++}^{2}$] and $\sigma_{L(5/2)}$ is not
differentiable on $-\mathbb{R}_{++}^{2}.$}
\end{example}

\begin{corollary}
\label{Saijo-FS}Assume that $F:\mathbb{R}_{+}^{p}\rightarrow\mathbb{R}_{+}$ is continuous.

\emph{(a)} If $F$ satisfies F.3c then $\sigma_{L(\gamma)}$ is differentiable
on $-\mathbb{R}_{++}^{p}$ for every $\gamma\in\Gamma_{F}$.

\emph{(b)} If $F$ satisfies F.2c and $\sigma_{L(\gamma)}$ is differentiable on
$-\mathbb{R}_{++}^{p}$ for every $\gamma\in\Gamma_{F}$, then $F$ satisfies F.3.
\end{corollary}

Proof. We know already that $L(\gamma)$ verifies condition (H) for every
$\gamma\in\Gamma_{F}.$

(a) Set $K:=\mathbb{R}_{+}^{p}$. Let $\gamma\in\Gamma_{F}$; since $F$ is
continuous and satisfies F.3c, we have that $A:=L(\gamma)$ is closed and
convex. By Proposition \ref{prop2} we have that $A$ verifies condition
$(\ref{r-sas})$. Applying Proposition \ref{prop-suf} we obtain that
$\sigma_{L(\gamma)}$ is differentiable on $-\mathbb{R}_{++}^{p}.$

(b) Take $\gamma\in\Gamma_{F}$. Since F.2c holds, by Proposition \ref{prop2}
we have that condition $(\ref{fp-ssc})$ holds for $A:=L(\gamma)$. Because
$\sigma_{A}$ is differentiable on $-\mathbb{R}_{++}^{p}$, using Theorem
\ref{Fact 14} we obtain that $A$ is convex. Hence $L(\gamma)$ is convex for
every $\gamma\in\Gamma_{F}$. Since $L(0)=\mathbb{R}_{+}^{p}$, $F$ is
quasi-concave, that is, F.3 holds. \hfill$\square$

\medskip

Having $F:\mathbb{R}_{+}^{p}\rightarrow\mathbb{R}_{+}$ a production function,
the cost function is defined for the price $x^{\ast}\in\mathbb{R}_{++}^{p}$
and the output $\gamma\in\mathbb{R}_{+}$ by
\[
c(x^{\ast},\gamma):=\inf\left\{  \left\langle x,x^{\ast}\right\rangle \mid
x\in L(\gamma)\right\}  =-\sigma_{L(\gamma)}(-x^{\ast}).
\]
So, the differentiability of $c(\cdot,\gamma)$ at $x^{\ast}\in\mathbb{R}%
_{++}^{p}$ (resp.\ on $\mathbb{R}_{++}^{p}$) is equivalent to the
differentiability of $\sigma_{L(\gamma)}$ at $-x^{\ast}\in-\mathbb{R}_{++}%
^{p}$ (resp.\ on $-\mathbb{R}_{++}^{p}$).

\begin{remark}
\label{rem4}\emph{Sakai \cite[Lem.\ 1 (2)]{Sak:73} obtained the
differentiability of the cost function for $F$ strictly concave and
satisfying F.1, F.2, F.4b; Saijo \cite{Sai:83} (see also Kim
\cite[Cor.\ 2]{Kim:93}) stated the same result for $F$ satisfying
F.2b, F.3b and F.4, but, as seen in Example \ref{ex-adsz}, this
result is not true; Fuchs-Selinger \cite[Thm.\ 2]{Fuc:97} obtained
the differentiability of the cost function for $F$ satisfying F.1,
F.2b, F.3b and F.4b, mentioning, using a figure, that the result of
Saijo is not true. Avriel \textit{et.~al.} \cite[Thm.\
4.8]{AvrDieSchZan:88} obtained the differentiability of the cost
function for $F$ satisfying F.1, F.2, F.3b and F.4b; our result in
Corollary \ref{Saijo-FS} is slightly more general.}
\end{remark}

\begin{example}
\label{ex28}\emph{Let the set $A\subset\mathbb{R}_{+}^{3}$
(considered in \cite{Zal:12}) be defined by
\begin{equation}
A:=\left\{  (x,y,z)\in\mathbb{R}_{+}^{3}\mid x+y\leq1,\ z\geq\frac{1}{2}%
\frac{(x+y-1)^{2}}{2-(x-y)^{2}}\right\}  \cup\left\{  (x,y,z)\in\mathbb{R}%
_{+}^{3}\mid x+y\geq1\right\}  . \label{d4}%
\end{equation}
The set $A$ is closed and convex, $0\notin A=A+K$, $\sigma_{A}$ is
differentiable on $\operatorname*{int}K^{-}$, but (\ref{rs2}) does
not hold (see \cite[Prop.\ 10]{Zal:12}). Moreover,
\begin{gather*}
\operatorname*{int}\nolimits_{\mathbb{R}_{+}^{p}}A=\left\{  (x,y,z)\in
\mathbb{R}_{+}^{3}\mid x+y\leq1,\ z>g(x,y)\right\}  \cup\left\{
(x,y,z)\in\mathbb{R}_{+}^{3}\mid x+y>1\right\}  ,\\
A+(\mathbb{R}_{+}^{p}\setminus\{0\})=\operatorname*{int}\nolimits_{\mathbb{R}%
_{+}^{p}}A,\quad\operatorname*{int}A=\mathbb{R}_{++}^{p}\cap
\operatorname*{int}\nolimits_{\mathbb{R}_{+}^{p}}A,\quad\mathbb{P}%
A=\mathbb{R}_{+}^{p}\setminus\{0\}.
\end{gather*}
It follows that $(\ref{fp-ssc})$ holds. Since $a:=(1,0,0)\in A$, $a^{\prime
}:=(0,1,0)\in A$ and $\tfrac{1}{2}a+\tfrac{1}{2}a^{\prime}\notin
\operatorname*{int}\nolimits_{\mathbb{R}_{+}^{p}}A$, we have that
$(\ref{r-sas})$ does not hold. Applying Theorem \ref{prop-fa} we have that the
corresponding function $F_{A}$ is a quasiconcave continuous production
function which verifies condition F.2c (even F.2b), but does not verify the
condition F.3c (that is $F_{A}$ is not strictly quasiconcave on $\mathbb{P}%
A$); consequently, $F_{A}$ is not strictly quasiconcave.}
\end{example}

\begin{remark}
\label{rem3}\emph{Corollary \ref{Saijo-FS}~(a) confirms the
\textquotedblleft only if\textquotedblright\ part of \cite[Thm.\
2]{Fuc:97}, while Corollary \ref{Saijo-FS}~(b), Example \ref{ex28}
and \cite[Prop.\ 10]{Zal:12} show that the \textquotedblleft
if\textquotedblright\ parts of \cite[Thm.\ 2]{Fuc:97} and of the
Equivalence Theorem in \cite{Sai:83} are not true; this is because
$F_{A}$ is continuous, satisfies F.2b and the corresponding cost
functions are differentiable on $-\mathbb{R}_{++}^{3}$, but $F_{A}$
is not strictly quasi-concave.}
\end{remark}

Note that \cite[Lem.\ 3]{Sai:83} follows from Proposition \ref{Fact 4}.
Indeed, because (USC) holds and $y\in\operatorname{Im}F$, $A:=L(y)\subset
\mathbb{R}_{+}^{p}$ is nonempty and closed; moreover, $\overline
{\operatorname*{conv}}\,A\subset\mathbb{R}_{+}^{p}$, and so $P:=(\overline
{\operatorname*{conv}}\,A)_{\infty}\subset\mathbb{R}_{+}^{p}$; hence
$\mathbb{R}_{+}^{p}\subset P^{+}$, and so $\mathbb{R}_{++}^{p}\subset
\operatorname*{int}P^{+}=\operatorname*{int}(\operatorname*{dom}g(\cdot,y)).$
In particular, we provided a new proof for Shephard's Lemma.

From Proposition \ref{Fact 4} we get also \cite[Thm.\ 1]{Fuc:95};
here $A$ is the set $R(x):=\{y\in X\mid y\succeq x\}$, $\succeq$
being a reflexive relation on the nonempty closed set
$X\subset\mathbb{R}_{+}^{p}$ such that $R(u)$ is closed for every
$u\in X$ (that is, $\succeq$ is upper semi\-continuous). Because any
nonempty closed subset of $\mathbb{R}_{+}^{p}$ can be represented as
$R(x)$ for some upper semi\-continuous reflexive relation $\succeq$
on $\mathbb{R}_{+}^{p}$, \cite[Thm.\ 2]{Fuc:95} is very close to
Proposition \ref{Fact 4} in the case the set $A$ is a subset of
$\mathbb{R}_{+}^{p}$ (because
$(\overline{\operatorname*{conv}}\,A)_{\infty
}\subset\mathbb{R}_{+}^{p}$).

\end{document}